\documentclass{amsart}
\usepackage{amsmath}
\usepackage{amsfonts}
\usepackage{enumerate}
\usepackage{amssymb}
\usepackage{amscd}
\usepackage[all]{xy}
\usepackage{bm}

\DeclareMathOperator{\gr}{gr} \DeclareMathOperator{\rk}{rk}
\DeclareMathOperator{\im}{Im} \DeclareMathOperator{\id}{id}
 
 \DeclareMathOperator{\Tor}{Tor}
\DeclareMathOperator{\Hom}{Hom} \DeclareMathOperator{\End}{End}
\DeclareMathOperator{\Ext}{Ext} \DeclareMathOperator{\Aut}{Aut}
\DeclareMathOperator{\Ann}{Ann} \DeclareMathOperator{\Ind}{Ind}
 \DeclareMathOperator{\Gal}{Gal}
\DeclareMathOperator{\GK}{GK} 
 \DeclareMathOperator{\reg}{reg}
\DeclareMathOperator{\Sym}{Sym} \DeclareMathOperator{\grm}{grmod}
\DeclareMathOperator{\Tot}{Tot}
\DeclareMathOperator{\degree}{degree} \DeclareMathOperator{\Ad}{Ad}
\DeclareMathOperator{\GL}{GL} \DeclareMathOperator{\Tr}{Tr}
\DeclareMathOperator{\cts}{cts}

\begin{document}
\theoremstyle{plain}
\newtheorem{MainThm}{Theorem}
\renewcommand{\theMainThm}{\Alph{MainThm}}
\newtheorem{MainCor}{Corollary}
\renewcommand{\theMainCor}{\Alph{MainCor}}
\newtheorem*{trm}{Theorem}
\newtheorem*{lem}{Lemma}
\newtheorem*{prop}{Proposition}
\newtheorem*{defn}{Definition}
\newtheorem*{thm}{Theorem}
\newtheorem*{example}{Example}
\newtheorem*{cor}{Corollary}
\newtheorem*{conj}{Conjecture}
\newtheorem*{hyp}{Hypothesis}
\newtheorem*{thrm}{Theorem}
\newtheorem*{quest}{Question}
\theoremstyle{remark}
\newtheorem*{rem}{Remark}
\newtheorem*{rems}{Remarks}
\newtheorem*{notn}{Notation}
\newcommand{\Fp}{\mathbb{F}_p}
\newcommand{\Fq}{\mathbb{F}_q}
\newcommand{\Zp}{\mathbb{Z}_p}
\newcommand{\Qp}{\mathbb{Q}_p}
\newcommand{\Kr}{\mathcal{K}}
\newcommand{\Rees}[1]{\widetilde{#1}}
\newcommand{\invlim}{\lim\limits_{\longleftarrow}}
\newcommand{\mgr}[1]{\mathcal{M}_{\gr}(#1)}
\newcommand{\grmod}[1]{\grm(#1)}

\title{$K_0$ and the dimension filtration for $p$-torsion Iwasawa modules}
\author{Konstantin Ardakov and Simon Wadsley}
\address{Ardakov:  Department of Pure Mathematics, University of Sheffield, Hicks Building, Hounsfield Road, Sheffield S3 7RH, UK}
\email{K.Ardakov@shef.ac.uk}
\address{Wadsley: DPMMS, University of Cambridge, Centre for Mathematical Sciences, Wilberforce Road, Cambridge CB3 0WB, UK}
\email{S.J.Wadsley@dpmms.cam.ac.uk}

\begin{abstract} Let $G$ be a compact $p$-adic analytic group. We study $K$-theoretic questions related to the representation theory of the completed group algebra $kG$ of $G$ with coefficients in a finite field $k$ of characteristic $p$. We show that if $M$ is a finitely generated $kG$-module whose dimension is smaller than the dimension of the centralizer of any $p$-regular element of $G$, then the Euler characteristic of $M$ is trivial. Writing $\mathcal{F}_i$ for the abelian category consisting of all finitely generated $kG$-modules of dimension at most $i$, we provide an upper bound for the rank of the natural map from the Grothendieck group of $\mathcal{F}_i$ to that of $\mathcal{F}_d$, where $d$ denotes the dimension of $G$. We show that this upper bound is attained in some special cases, but is not attained in general.
\end{abstract}

\thanks{This work was started when the first author was the Sir Robert and Lady Clayton Junior Research Fellow at Christ's College, Cambridge}

\maketitle
\let\le=\leqslant  \let\leq=\leqslant
\let\ge=\geqslant  \let\geq=\geqslant

\section{Introduction}

\subsection{Iwasawa algebras}
\label{start} In this paper we study certain aspects of the representation theory of Iwasawa algebras. These are the completed group algebras
\[ \Lambda_G := \invlim \Zp[G/U], \]
where $\Zp$ denotes the ring of $p-$adic integers, $G$ is a compact $p-$adic analytic group, and the inverse limit is taken over the open normal subgroups $U$ of $G$. Closely related is the epimorphic image $\Omega_G$ of $\Lambda_G$,
\[\Omega_G := \invlim \Fp[G/U], \]
where $\Fp$ is the field of $p$ elements.

This paper is a continuation of our earlier work \cite{AW}, in which
we investigated the relationship between the notion of
characteristic element of $\Lambda_G$-modules defined in
\cite{CFKSV} and the Euler characteristic of $\Omega_G$-modules. We
focus exclusively on the $p$-torsion $\Lambda_G$-modules, that is,
those killed by a power of $p$. Because we are interested in
$K$-theoretic questions, we only need to consider those modules
actually killed by $p$, or equivalently, the $\Omega_G$-modules.

In this introduction we assume that the group $G$ has no elements of
order $p$, although all of our results hold for arbitrary $G$ with slightly more involved formulations.

\subsection{The dimension filtration}
\label{IntroDim} The category $\mathcal{M}(\Omega_G)$ of finitely
generated $\Omega_G$-modules has a canonical dimension function $M\mapsto
d(M)$ defined on it which provides a filtration of $\mathcal{M}(\Omega_G)$
by admissible subcategories $\mathcal{F}_i$ whose objects are those
modules of dimension at most $i$. If $d$ denotes the dimension of
$G$ then $\mathcal{F}_d$ is just $\mathcal{M}(\Omega_G)$ and
$\mathcal{F}_{d-1}$ is the full subcategory of torsion modules.

Perhaps the central result of \cite{AW} was to classify those
$p$-adic analytic groups $G$ for which the Euler characteristic of
every finitely generated torsion $\Omega_G$-module is trivial. Following
the proof of \cite[Theorem 8.2]{AW} we see that the Euler
characteristic of every module $M \in \mathcal{F}_i$ is trivial if
and only if the natural map from $K_0(\mathcal{F}_i)$ to
$K_0(\mathcal{F}_d)$ is the zero map. This raises the following

\begin{quest} When is the natural map $\alpha_i : K_0(\mathcal{F}_i) \to K_0(\mathcal{F}_d)$ the zero map?
\end{quest}

Let $\Delta^+$ denote the largest finite normal subgroup of $G$. In
\cite{AW} we answered the above question in the case $i=d-1$: the
map $\alpha_{d-1}$ is zero precisely when $G$ is
\emph{$p$-nilpotent}, that is, when $G/\Delta^+$ is a pro-$p$ group.
There is a way of rephrasing this condition in terms of $G_{\reg}$,
the set of elements of $G$ of finite order: $G$ is $p$-nilpotent if
and only if the centralizer of every element $g \in G_{\reg}$ is an
open subgroup of $G$.

\subsection{Serre's work}
\label{SerreWork} This question was answered by Serre in
\cite{Serre3} in the case $i=0$, although he did not use our
language. In this case $\mathcal{F}_i$ consists of modules that are
finite dimensional as $\Fp$-vector spaces. He produced a formula which
relates the Euler characteristic $\chi(G,M)$ of a module $M \in
\mathcal{F}_0$ with its Brauer character $\varphi_M$:
\[\log_p \chi(G,M) = \int_G \varphi_M(g) \det(1 - \Ad(g^{-1})) dg.\]
Here $\Ad : G \to \GL(\mathcal{L}(G))$ is the adjoint representation
of $G$. As a consequence, Serre proved that the Euler characteristic
of every module $M \in \mathcal{F}_0$ is trivial precisely when the
centralizer of every element $g \in G_{\reg}$ is infinite.

\subsection{Trivial Euler characteristics}

The results of Serre and \cite{AW} mentioned above suggest that
there might be a connection between the answer to Question
\ref{IntroDim} and the dimensions of the centralizers of the
elements of $G_{\reg}$. Indeed, one might wonder whether $\alpha_i$
is zero if and only if the $\dim C_G(g) > i$ for all elements $g \in
G_{\reg}$. Whilst this latter statement turns out to be not quite
right --- see (\ref{CEx}) --- such a connection indeed exists.

\begin{MainThm} \label{TrivialEul} Let $M$ be a finitely generated $\Omega_G$-module such that $d(M) < \dim C_G(g)$ for all $g \in G_{\reg}$. Then $\chi(G,M) = 1$.
\end{MainThm}
The proof is given in (\ref{TrivEulCh}).

\subsection{A related question}
\label{RankQ}
In the remainder of the paper, we address the following
\begin{quest} What is the rank of the natural map $\alpha_i : K_0(\mathcal{F}_i) \to K_0(\mathcal{F}_d)$?
\end{quest}

The group $G$ acts on $G_{\reg}$ by conjugation and this action commutes with the action of a certain Galois group $\mathcal{G}$ on $G_{\reg}$, which essentially acts by raising elements to powers of $p$; see (\ref{preg}) for details. Let $S_i=\{g\in G_{\reg} : \dim C_G(g) \leq i\}$  --- this is a union of $G \times \mathcal{G}$-orbits. For example, $S_d = G_{\reg}$ and $G_{\reg}  - S_{d-1}$ is the set of all elements of $G$ which have finite order and lie in a finite conjugacy class; in fact, $G_{\reg} - S_{d-1}$ is just the largest finite normal subgroup $\Delta^+$ of $G$ mentioned in (\ref{IntroDim}).

Our next result provides an upper bound for the rank of $\alpha_i$:

\begin{MainThm} \label{UpBnd} $\rk \alpha_i$ is bounded above by the number of $G\times\mathcal{G}$-orbits on $S_i$.
\end{MainThm}
The proof is given in (\ref{UpBndPf}).

\subsection{Some special cases}
\label{MainQuest}
We next show that the rank of $\alpha_i$ attains the upper bound given in Theorem \ref{UpBnd} in some special cases.

\begin{MainThm} The rank of $\alpha_i$ equals the number of $G\times\mathcal{G}$-orbits on $S_i$ if either
\begin{enumerate}[{(}a{)}]
\item $i=d$, or $i=d-1$ or $i=0$, or if
\item $G$ is virtually abelian.
\end{enumerate}
\end{MainThm}
This follows by combining Propositions
\ref{KoFd},\ref{torsion},\ref{rkalpha0} and Theorem \ref{VirtAb}.

Finally, in (\ref{CEx}) we give an example to show that the upper bound of Theorem \ref{UpBnd} is \emph{not} always attained. Questions \ref{IntroDim} and \ref{RankQ} remain open in general.

\section{Generalities}

Throughout this paper, $k$ will denote a fixed finite field of characteristic $p$ and order $q$. Modules will be right modules, unless explicitly stated otherwise. We will conform with the notation of \cite{AW}, with the exception of \cite[\S 12]{AW}.

\subsection{Grothendieck groups}\label{GrotGps} Let $\mathcal{A}$ be a small abelian category. A full additive subcategory $\mathcal{B}$ of $\mathcal{A}$ is \emph{admissible} if whenever $0 \to M' \to M \to M'' \to 0$ is a short exact sequence in $\mathcal{A}$ such that $M$ and $M''$ belong to $\mathcal{B}$, then $M'$ also belongs to $\mathcal{B}$ \cite[12.4.2]{MCR}.

The \emph{Grothendieck group} $K_0(\mathcal{B})$ of $\mathcal{B}$ is the abelian group with generators $[M]$ where $M$ runs over all the objects of $\mathcal{B}$ and relations $[M] = [M'] + [M'']$ for any short exact sequence $0 \to M' \to M \to M'' \to 0$ in $\mathcal{A}$ \cite[12.4.3]{MCR}.

We will frequently be dealing with vector space versions of these groups. To simplify the notation later on, whenever $F$ is a field we will write
\[FK_0(\mathcal{B}) := F \otimes_{\mathbb{Z}} K_0(\mathcal{B}).\]

If $A$ is a ring, then $\mathcal{P}(A)$, the category of all finitely generated projective modules is an admissible subcategory of $\mathcal{M}(A)$, the category of all finitely generated $A$-modules. The \emph{Grothendieck groups} of $A$ are defined as follows:
\begin{itemize}
\item $K_0(A) := K_0(\mathcal{P}(A))$, and
\item $\mathcal{G}_0(A) := K_0(\mathcal{M}(A))$.
\end{itemize}
We also set $FK_0(A) := FK_0(\mathcal{P}(A))$ and $F\mathcal{G}_0(A) := FK_0(\mathcal{M}(A))$.

\subsection{Homology and Euler characteristics}
\label{HomEul}
Let $\mathcal{A}$ be an abelian category, let $\Gamma$ be an abelian group and let $\psi$ be an additive function from the objects of $\mathcal{A}$ to $\Gamma$. This means that for every object $A \in \mathcal{A}$ there exists an element $\psi(A) \in \Gamma$ such that
\begin{itemize}
\item $\psi(B) = \psi(A) + \psi(C)$ whenever $0\to A\to B\to C\to 0$ is a short exact sequence in $\mathcal{A}$.
\end{itemize}
In this context, if $C_\ast = \cdots \to C_n \to \cdots \to C_0 \to \cdots$ is a bounded complex in $\mathcal{A}$ (that is, $C_i = 0$ for sufficiently large $|i|$), we define the \emph{Euler characteristic} of $C_\ast$ (with respect to $\psi$) to be
\[\psi(C_\ast) := \sum_{i\in\mathbb{Z}}(-1)^i \psi(C_i) \in \Gamma.\]
Note that if we think of any object $A \in \mathcal{A}$ as a complex $A_\ast$ concentrated in degree zero, then $\psi(A_\ast) = \psi(A)$. The following is well-known.

\begin{lem} Let $C_\ast$ be a bounded complex in $\mathcal{A}$. Then the Euler characteristic of $C_\ast$ equals the Euler characteristic of the homology complex $H_\ast(C_\ast)$ with zero differentials:
\[\psi(C_\ast) = \psi(H_\ast(C_\ast)).\]
\end{lem}

Thus $\psi$ extends to a well-defined function from the objects of the derived category $\mathbf{D}(\mathcal{A})$ of $\mathcal{A}$ (if this exists) to $\Gamma$.

\subsection{Spectral sequences and Euler characteristics}
\label{SpecEul}
We will require a version of Lemma \ref{HomEul} for spectral sequences. Let $E = E^r_{ij}$ be a homology spectral sequence in $\mathcal{A}$ starting with $E^a$ \cite[\S 5.2]{Wei}. We say that $E$ is \emph{totally bounded} if $E^a_{ij}$ is zero for all sufficiently large $|i|$ or $|j|$. This is equivalent to insisting that $E$ is bounded in the sense of \cite[5.2.5]{Wei} and that there are only finitely many nonzero diagonals on each page. It is clear that $E$ converges.

We define the \emph{$r^{th}$-total complex} $\Tot(E^r)_\ast$ of $E$ to be the complex in $\mathcal{A}$ whose $n^{th}$ term is
\[\Tot(E^r)_n = \bigoplus_{i+j=n}E^r_{ij}\]
and whose $n^{th}$ differential is $\bigoplus_{i \in \mathbb{Z}} d^r_{i,n-i}$. Because $E$ is totally bounded, $\Tot(E^r)_n \in \mathcal{A}$ for all $r\geq a$ and $n\in \mathbb{Z}$ and each complex $\Tot(E^r)$ is bounded.

The next result is folklore --- see for example \cite[Example 6, p.15]{McC} for a cohomological version and the proof of \cite[Theorem 9.8]{Lang} for a similar formulation --- but we give a proof for the convenience of the reader.

\begin{prop} Let $E$ be a totally bounded homology spectral sequence in $\mathcal{A}$ starting with $E^a$. Then $\psi(\Tot(E^a)_\ast) = \psi(\Tot(E^\infty)_\ast)$.
\end{prop}
\begin{proof}
From the definition of a spectral sequence we see that
\[H_n(\Tot(E^r)_\ast) = \Tot(E^{r+1})_n\]
for all $r \geq a$ and $n\in\mathbb{Z}$. Applying Lemma \ref{HomEul} repeatedly gives
\[\psi(\Tot(E^a)_\ast) = \psi(\Tot(E^{a+1})_\ast) = \cdots = \psi(\Tot(E^\infty)_\ast)\]
as required.
\end{proof}

\subsection{Completed group algebras} Let $G$ be a profinite group. We will write $\mathbf{U}_G$
(respectively, $\mathbf{U}_{G,p}$) for the set of all open normal
(respectively, open normal pro-$p$) subgroups $U$ of $G$.

Define the \emph{completed group algebra} of $G$ by the formula
\[kG := k[[G]] := \lim\limits_{\stackrel{\longleftarrow}{U\in\mathbf{U}_G}}k[G/U].\]
As each group algebra $k[G/U]$ is finite, this is a compact
topological $k$-algebra which ``controls" the continuous
$k$-representations of $G$ in the following sense. Whenever $V$ is a
compact topological $k$-vector space and $\rho : G \to
\Aut_{\cts}(V)$ is a continuous representation of $G$ then $\rho$
extends to a unique continuous homomorphism of topological
$k$-algebras $\rho : kG \to \End_{\cts}(V)$, and $V$ becomes a
compact topological (left) $kG$-module. Conversely, any compact
topological $kG$-module $V$ gives rise to a continuous
$k$-representation of $G$.

When $G$ is a compact $p$-adic analytic group, $kG$ is sometimes
called an \emph{Iwasawa algebra} --- we refer the reader to
\cite{AB2} for more details. We note that in this case $kG$ is
always Noetherian and has finite global dimension when $G$ has no
elements of order $p$ --- facts which we will sometimes use without
further mention.

\section{Brauer characters}
Throughout this section, $G$ denotes a profinite group which is virtually pro-$p$.

\label{BrChSec}
\subsection{$\bm{p}$-regular elements and the Galois action}
\label{preg}
An element $g \in G$ is said to be \emph{$p$-regular} if its order, in the profinite sense, is coprime to $p$. We will denote the set of all $p$-regular elements of $G$ by $G_{\reg}$ --- this is a union of conjugacy classes in $G$. Because $G$ is virtually pro-$p$, any $p$-regular element has finite order; moreover, if $G$ has no elements of order $p$ then $G_{\reg}$ is just the set of elements of finite order in $G$, so this definition extends the one given in (\ref{IntroDim}).

Let $m$ denote the $p'$-part of $|G|$ in the profinite sense; equivalently, $m$ is the index of a Sylow pro-$p$ subgroup of $G$. As we are assuming that $G$ is virtually pro-$p$, $m$ is finite. Let $k' = k(\omega)$, where $\omega$ is a primitive $m$-th root of unity over $k$ and let $\mathcal{G}_k$ be the Galois group $\Gal(k(\omega)/k)$. If $\sigma \in \mathcal{G}_k$, then $\sigma(\omega) = \omega^{t_{\sigma}}$ for some $t_{\sigma} \in (\mathbb{Z}/m\mathbb{Z})^\times$. This gives an injection $\sigma \mapsto t_{\sigma}$ of $\mathcal{G}_k$ into $(\mathbb{Z}/m\mathbb{Z})^\times$.

We can now define a left permutation action of $\mathcal{G}_k$ on
$G_{\reg}$ by setting $\sigma.g = g^{t_{\sigma}}$; this makes sense
because $t_{\sigma}$ is coprime to the order of any element $g\in
G_{\reg}$ by construction. This action commutes with any
automorphism of $G$, so $\mathcal{G}_k$ permutes the $p$-regular
conjugacy classes of $G$. These constructions give a continuous
action $G\times \mathcal{G}_k$ on $G_{\reg}$; note that
$G\times\mathcal{G}_k$ is also virtually pro-$p$ because
$\mathcal{G}_k$ is a finite group.

\subsection{Locally constant functions}
Let $X$ be a compact totally disconnected space. For any commutative ring $A$ we let $C(X;A)$ denote the ring of all locally constant functions from $X$ to $A$. If $G$ acts on $X$ continuously on the left then it acts on $C(X;A)$ on the right as follows:
\[(f.g)(x) = f(g.x)\quad\mbox{for all}\quad f \in C(X;A), g \in G, x \in X.\]
We will identify the subring of invariants $C(X;A)^G$ with $C(G \backslash X;A)$, where $G \backslash X$ denotes the set of $G$-orbits in $X$.

\subsection{Brauer characters}
\label{BrCh}
Our treatment is closely follows Serre \cite[\S 2.1, \S 3.3]{Serre3}.

Fix a finite unramified extension $K$ of $\Qp$ with residue field $k$. Let $F = K(\tilde{\omega})$, where $\tilde{\omega}$ is a primitive $m$-th root of $1$. Then the ring of integers of $F$ is $\mathcal{O}' = \mathcal{O}[\tilde{\omega}]$ where $\mathcal{O}$ is the ring of integers of $K$. Reduction modulo $p$ gives an isomorphism of the residue field of $F$ with $k'$ and we may assume that $\tilde{\omega}$ maps to $\omega$ under this isomorphism. For each $m$-th root of unity $\xi \in k'$ there is a unique $m$-th root of unity $\tilde{\xi} \in F$ such that $\tilde{\xi}$ maps to $\xi$ modulo $p$. This gives us an isomorphism $\widetilde{\hspace{0.1cm}\cdot\hspace{0.1cm}}: \langle \omega \rangle \to \langle \tilde{\omega}\rangle $ between the two cyclic groups.

Let $\mathcal{F}_0 = \mathcal{F}_0(G)$ be the abelian category of all topological $kG$-modules which are finite dimensional over $k$. If $A \in \mathcal{F}_0$ and $g \in G_{\reg}$, the eigenvalues of the action of $g$ on $A$ are powers of $\omega$ --- say $\xi_1,\ldots,\xi_d$ (always counted with multiplicity so that $\dim A = d$). Define
\[\varphi_A(g) = \sum_{i=1}^d \widetilde{\xi_i} \in F.\]
The function $\varphi_A : G_{\reg} \to F$ is called the $\emph{Brauer character}$ of $A$. It has the following properties:

\begin{lem} Let $A,B,C \in \mathcal{F}_0$.
\begin{enumerate}[{(}i{)}]
\item $\varphi_A$ is a locally constant $G\times\mathcal{G}_k$-invariant $F$-valued function on $G_{\reg}$:
\[\varphi_A \in C(G_{\reg};F)^{G\times\mathcal{G}_k}.\]
\item $\varphi_B = \varphi_A + \varphi_C$ whenever $0\to A \to B \to C \to 0$ is a short exact sequence.
\item $\varphi_{A\otimes_k B} = \varphi_A \cdot \varphi_B$, where $G$ acts diagonally on $A\otimes_k B$.
\item Let $A' = A\otimes_k k'$ be the $k'H$-module obtained from $A$ by extension of scalars. Then $\varphi_{A'} = \varphi_A$.
\end{enumerate}
\end{lem}
\begin{proof} As $G$ acts continuously on the finite dimensional vector space $A$, some $U\in\mathbf{U}_G$ acts trivially. It follows that for any $g\in G$, $\varphi_A$ is constant on the open neighbourhood $gU$ of $g$, so $\varphi_A$ is a locally constant function. For the remaining assertions, we may assume that $G$ is finite. In this case, the result is well-known --- see for example \cite[Volume I, \S 17A, \S 21B]{CurRei} or \cite[\S 18]{Serre2}.
\end{proof}

\subsection{Berman--Witt Theorem}
\label{BermanWitt} Lemma \ref{BrCh} shows that there is an $F$-linear map
\[\varphi : FK_0(\mathcal{F}_0) \to C(G_{\reg};F)^{G\times\mathcal{G}_k}\]
given by $\varphi(\lambda \otimes [A]) = \lambda\varphi_{A}$ for all $\lambda \in F$ and $A \in \mathcal{F}_0$ --- see (\ref{GrotGps}) for the notation.

\begin{thm}
$\varphi$ is an isomorphism.
\end{thm}

This is a generalization of the well-known Berman--Witt Theorem
\cite[Volume I, Theorem 21.25]{CurRei} in the case when $G$ is
finite. When $k = \Fp$ and $G$ is finite, a short proof can also be
found in \cite[\S 2.3]{Serre3}. Our proof is given below in
(\ref{PfBmWt}).

Until the end of $\S \ref{BrChSec}$, we fix $U\in\mathbf{U}_{G,p}$ and write $\overline{G}$ for the quotient group $G/U$.

\subsection{Lemma}
\label{Piuonto}
Let $\pi_U : G \twoheadrightarrow \overline{G}$ be the natural surjection. Then
\[\pi_U(G_{\reg}) = \overline{G}_{\reg}.\]
\begin{proof}
It is enough to show that $\pi_U(G_{\reg}) \supseteq \overline{G}_{\reg}$.  First suppose that $G$ is a finite group so that $U$ is a normal $p$-subgroup. For any $x \in G$ we can find unique commuting $s \in G_{\reg}$ and $u \in G$ such that $x = su$ and $u$ has order a power of $p$. Now if $xU \in \overline{G}$ is $p$-regular then raising $x$ to an appropriate sufficiently large power of $p$ doesn't change $xU$ and has the effect of making $x = s$ $p$-regular.

Now suppose that $G$ is arbitrary. Serre \cite[\S 1.1]{Serre3} has observed that $G_{\reg}$ is a compact subset of $G$ which can be identified with the inverse limit of the various $(G/W)_{\reg}$ as $W$ runs over $\mathbf{U}_G$. Because $U$ is open in $G$, we may assume that all the $W$'s are contained in $U$. Now the result follows from the first part.
\end{proof}

\subsection{Proposition}\label{HallBij} The map $\pi_U$ induces a bijection
\[\pi_U : (G\times\mathcal{G}_k) \backslash G_{\reg} \to (\overline{G}\times\mathcal{G}_k) \backslash \overline{G}_{\reg}.\]
\begin{proof}
In view of Lemma \ref{Piuonto}, it is sufficient to prove that this
map is injective. So let $x,y\in G_{\reg}$ be such that $xU$ and
$yU$ lie in the same $\overline{G} \times \mathcal{G}_k$-orbit. By
replacing $y$ by a $G\times\mathcal{G}_k$-conjugate, we may assume
that actually $xU = yU$. As $G$ is virtually pro-$p$, it will now be
sufficient to show that $xW$ and $yW$ are conjugate in $G/W$ for any
$W\in\mathbf{U}_{G,p}$ contained in $U$. Without loss of generality,
we can assume that $G$ is finite and that $W=1$.

Now as $x$ is $p$-regular and $U$ is a $p$-group, $\langle x\rangle \cap U = 1$. Similarly $\langle y\rangle \cap U=1$, so $\langle x\rangle \cong \langle xU\rangle = \langle yU\rangle \cong \langle y\rangle$.

Consider the finite solvable group $H := \langle x\rangle U$. As $xU
= yU$, $y$ lies in $H$ and $U$ is the unique Sylow $p$-subgroup of
$H$. It follows that $\langle x\rangle$ and $\langle y\rangle$ are
\emph{Hall} $p'$-\emph{subgroups} of $H$ and as such are conjugate
in $H$ \cite{Hall}. Hence there exists $h \in H$ such that $x^h :=
h^{-1}xh = y^a$ for some $a\geq 1$. Now as $H/U$ is abelian, $xU =
x^h U$. Hence $yU = y^a U$, but $y^i \mapsto y^iU$ is an isomorphism
so $y = y^a = x^h$ as required.
\end{proof}

We would like to thank Jan Saxl for providing this proof.

\begin{cor} The $G\times \mathcal{G}_k$-orbits in $G_{\reg}$ are closed and open in $G_{\reg}$.
\end{cor}
\begin{proof}
Because $G\times\mathcal{G}_k$ is a profinite group acting continuously on the Hausdorff space $G_{\reg}$, the orbits are closed. But they are disjoint and finite in number by the Proposition, so they must also be open.
\end{proof}

\subsection{D\'evissage} \label{Dev}
Any finitely generated $k\overline{G}$-module is finite dimensional over $k$ since $\overline{G}$ is finite, so we have a natural inclusion $\mathcal{M}(k\overline{G}) \subset \mathcal{F}_0$ of abelian categories. Let $\lambda_U : \mathcal{G}_0(k\overline{G}) \to K_0(\mathcal{F}_0)$ be the map induced on Grothendieck groups.

\begin{lem} $\lambda_U$ is an isomorphism.
\end{lem}
\begin{proof}
Let $w_U := (U-1)kG$ be the kernel of the natural map $kG \twoheadrightarrow k\overline{G}$ and let $A \in \mathcal{F}_0$. As in the proof of Lemma \ref{BrCh} we can find $W\in\mathbf{U}_G$ which acts trivially on $A$, which we may assume to be contained in $U$. Now $A$ is a $k[G/W]$-module and the image of $w_U$ in $k[G/W]$ is nilpotent because $U/W$ is a normal $p$-subgroup of the finite group $G/W$ \cite[Lemma 3.1.6]{Pass2}. Thus $Aw_U^t = 0$ for some $t\geq 0$, so $A$ has a finite filtration
\[0 = Aw_U^t \subset Aw_U^{t-1} \subset\cdots\subset Aw_U\subset A\]
where each factor is a $k\overline{G}$-module. Hence $\lambda_U$ is an isomorphism by d\'evissage --- see, for example, \cite[Theorem 12.4.7]{MCR}.
\end{proof}
\subsection{Proof of Theorem \ref{BermanWitt}}
\label{PfBmWt}
Consider the commutative diagram
\[
\xymatrix{
F\mathcal{G}_0(k\overline{G}) \ar[d]_{\lambda_U} \ar[r]^{\varphi} & C(\overline{G}_{\reg};F)^{ \overline{G}\times\mathcal{G}_k} \ar[d]^{\pi_U^\ast}\\
FK_0(\mathcal{F}_0) \ar[r]^{\varphi} & C(G_{\reg};F)^{ G\times\mathcal{G}_k},\\
}
\]
where the map $\pi_U^\ast$ is defined by the formula $\pi_U^\ast(f)(g) = f(gU)$. The top horizontal map $\varphi$ is an isomorphism by the usual Berman--Witt Theorem \cite[Volume I, Theorem 21.25]{CurRei}, $\pi_U^\ast$ is an isomorphism as a consequence of Proposition \ref{HallBij} and $\lambda_U$ is an isomorphism by Lemma \ref{Dev}. Hence the bottom horizontal map $\varphi$ is also an isomorphism, as required. \qed

\section{Modules over Iwasawa algebras}
\label{ModsIwa}
\subsection{Compact $p$-adic analytic groups}
From now on, $G$ will denote an arbitrary compact $p$-adic analytic group. We will write $d := \dim G$ for the dimension of $G$. By the celebrated result of Lazard \cite[Corollary 8.34]{DDMS}, $G$ has an open normal \emph{uniform} pro-$p$ subgroup, so $G$ is in particular virtually pro-$p$ and we can apply the theory developed in \S \ref{BrChSec}. We fix such a  subgroup $N$ in what follows, and write $\overline{G}$ for the quotient group $G/N$.

\subsection{The base change map}
\label{BaseChange}
We begin by studying $\mathcal{G}_0(kG)$ in detail. First, a preliminary
result.

\begin{lem} The projective dimension of the left $kG$-module $k\overline{G}$ is at most $d$.
\end{lem}
\begin{proof}
It is well known that $kN$ has global dimension $\dim N = d$ \cite[Theorem 4.1]{Bru}. Now $k\overline{G} \cong kG\otimes_{kN} k$ as a left $kG$-module and $kG$ is a free right $kN$-module of finite rank. The result follows.
\end{proof}

Hence the right $k\overline{G}$-modules $\Tor_j^{kG}(M,k\overline{G})$ are zero for all $j>d$ and all $M\in\mathcal{M}(kG)$, so we can define an element $\theta_N[M] \in \mathcal{G}_0(k\overline{G})$ by the formula
\[\theta_N[M]:=\sum_{j=0}^\infty (-1)^j[\Tor_j^{kG}(M,k\overline{G})].\]
The long exact sequence for $\Tor$ shows that $M \mapsto \theta_N[M]$ is an additive function on the objects of $\mathcal{M}(kG)$, so we have a base change map \cite[p. 454]{Bass} on $\mathcal{G}$-theory
\[\theta_N : \mathcal{G}_0(kG) \to \mathcal{G}_0(k\overline{G})\]
that will be one of our tools for studying $\mathcal{G}_0(kG)$.

\subsection{Graded Brauer characters}
\label{GrBr}
Let $\grmod{kG}$ denote the category of all $kG$-modules $M$ which admit a direct sum decomposition
\[M = \bigoplus_{n\in\mathbb{Z}} M_n\]
into a $kG$-submodules such that $M_n$ is finite dimensional for all $n\in \mathbb{Z}$ and zero for all sufficiently small $n$, thought of as a full subcategory of the abelian category of all $\mathbb{Z}$-graded $kG$-modules and graded maps of degree zero.

\begin{defn} The \emph{graded Brauer character} of $M \in \grmod{kG}$ is the function
\[\zeta_M : G_{\reg} \to F[[t,t^{-1}]]\]
defined by the formula
\[\zeta_M(g) = \sum_{n\in\mathbb{Z}} \varphi_{M_n}(g) t^n.\]
\end{defn}
This definition extends the notion of Brauer character presented in (\ref{BrCh}) if we think of any finite dimensional $kG$-module $M$ as a graded module concentrated degree zero.

Now let $M \in \mathcal{M}(kG)$ and let $w_N := (N-1)kG$ be the
kernel of the natural map $kG \twoheadrightarrow k\overline{G}$. As
 $kG$ is Noetherian, $w_N^n$ is a finitely generated right ideal in
$kG$ for all $n\geq 0$, so the modules $Mw_N^n/Mw_N^{n+1}$ are
finite dimensional over $k$ for all $n$. Hence the associated graded
module
\[\gr M := \bigoplus_{n=0}^\infty \frac{Mw_N^n}{Mw_N^{n+1}}\]
lies in $\grmod{kG}$, and as such has a graded Brauer character $\zeta_{\gr M}$.

We will see in (\ref{BrRat}) that $\zeta_{\gr M}(g)$ is actually a
rational function in $t$ for each $g \in G_{\reg}$.

\subsection{The adjoint representation}
\label{AdRep}
Recall \cite[\S 4.3]{DDMS} that there is an \emph{additive structure} $(N,+)$ on our fixed uniform subgroup $N$. In this way $N$ becomes a free $\Zp$-module of rank $d$ so $\mathcal{L}(G) := \Qp \otimes_{\Zp} N$ is a $\Qp$-vector space of dimension $d$. There is a way of turning $\mathcal{L}(G)$ into a Lie algebra over $\Qp$ \cite[\S 9.5]{DDMS}, but we will not need this.

The conjugation action of $G$ on $N$ respects the additive structure on $N$ and gives rise to the \emph{adjoint representation}
\[\Ad : G \to \GL(\mathcal{L}(G))\]
given by $\Ad(g)(n) = gng^{-1}$ for all $g \in G$ and $n \in N$. We define a function $\Psi:G \to \Qp[t]$ by setting
\[\Psi(g) := \det(1 - \Ad(g^{-1})t)\]
for all $g \in G$. As $\Psi(g)\cdot \det\Ad(g)$ is the characteristic polynomial of $\Ad(g)$, we can think of $\Psi(g)$ as a polynomial in $F[t]$ of degree $d$.

\subsection{The key result}
\label{MainResult}
By Theorem \ref{BermanWitt} and Lemma \ref{Dev}, the composite map $\varphi\circ\lambda_N : F\mathcal{G}_0(k\overline{G}) \to C(G_{\reg};F)^{G\times \mathcal{G}_k}$ is an isomorphism. We therefore do not lose any information when studying $\theta_N$ by postcomposing it with this isomorphism. Our main technical result reads as follows:

\begin{thm} Let $\rho_N$ be the composite map
\[\rho_N := \varphi \circ \lambda_N \circ \theta_N : F\mathcal{G}_0(kG) \to C(G_{\reg};F)^{G\times \mathcal{G}_k}.\]
Then for any $g \in G_{\reg}$ and $M \in \mathcal{F}_d$, the number
\[\rho_N[M](g) = \sum_{j=0}^d (-1)^j \varphi_{\Tor^{kG}_j(M,k\overline{G})}(g) \in F\]
equals the value at $t=1$ of the rational function
\[\zeta_{\gr M}(g) \cdot \Psi(g) \in F(t).\]
\end{thm}

We now begin preparing for the proof, which is given in
$\S\ref{PhiZeroProof}$.

\section{Graded modules for $\Sym(V) \# H$}
\label{GrMod}
\subsection{Notation} Let $V$ be a finite dimensional $k$-vector space and let $H$ be a finite group acting on $V$ by $k$-linear automorphisms on the right. We will write $v^h$ for the image of $v \in V$ under the action of $h\in H$. This action extends naturally to an action of $H$ on the symmetric algebra $\Sym(V)$ by $k$-algebra automorphisms. Let
\[R := \Sym(V) \# H\]
denote the skew group ring \cite[1.5.4]{MCR}: by definition, $R$ is a free right $\Sym(V)$-module with basis $H$, with multiplication given by the formula
\[ (hr)(gs) = (hg)(r^gs)\]
for all $g,h \in H$ and $r,s \in \Sym(V)$. Thus $R$ is isomorphic to $kH\otimes_k \Sym(V) $ as a $k$-vector space and setting $R_n := kH\otimes_k\Sym^nV$ turns $R = \bigoplus_{n=0}^\infty R_n$ into a graded $k$-algebra.

\subsection{Dimensions}
\label{Dims}
If $S$ is a positively graded $k$-algebra, let $\mgr{S}$ denote the category of all finitely generated $\mathbb{Z}$-graded right $S$-modules and graded maps of degree zero. Since $R = \Sym(V)\#H$ is a finitely generated $\Sym(V)$-module, we have an inclusion $\mgr{R} \subset \mgr{\Sym(V)}$ of abelian categories.

Following \cite[\S 11]{AM} we define the \emph{dimension} $d(M)$ of a module $M\in\mgr{\Sym(V)}$ to be the order of the pole of the \emph{Poincar\'e series} $P_M(t)$ of $M$ at $t=1$, where
\[P_M(t) = \sum_{n\in\mathbb{Z}} (\dim_k M_n) t^n \in \mathbb{Z}[[t,t^{-1}]].\]
In fact \cite[Theorem 11.1]{AM}, there exists a polynomial $u(t) \in \mathbb{Z}[t,t^{-1}]$ such that $u(1) \neq 0$ and
\[P_M(t) = \frac{u(t)}{(1-t)^{d(M)}}.\]
Note that we have to allow Laurent polynomials and power series because our modules are $\mathbb{Z}$-graded. It is well-known that $d(M)$ also equals the Krull dimension $\Kr(M)$ of $M$ in the sense of \cite[\S 6.2]{MCR} and the Gelfand-Kirillov dimension $\GK(M)$ of $M$ \cite[\S 8]{MCR}.

\subsection{Properties of graded Brauer characters}
\label{GrBrProps}

Recall the definition of graded Brauer characters given in
(\ref{GrBr}). As the finite group $H$ is in particular compact
$p$-adic analytic, we may speak of the category $\grmod{kH}$. The
following result is a straightforward application of Lemma
\ref{BrCh} and shows that graded Brauer characters behave well with
respect to basic algebraic constructions.

\begin{lem} Let $A,B,C \in \grmod{kH}$.
\begin{enumerate}[{(}i{)}]
\item If $0 \to A \to B \to C \to 0$ is an exact sequence in $\grmod{kH}$ then
\[\zeta_B = \zeta_A + \zeta_C.\]
\item Define $A \otimes_k B \in \grmod{kH}$ by letting $H$ act diagonally and giving $A\otimes_k B$ the tensor product gradation
\[(A \otimes_k B)_n = \bigoplus_{i+j=n} A_i \otimes_k B_j.\]
Then $\zeta_{A\otimes_k B} = \zeta_A \cdot \zeta_B.$
\item For each $m\in\mathbb{Z}$ define the \emph{shifted module} $A[m] \in \grmod{kH}$ by setting $A[m]_n = A_{m+n}$ for all $n\in\mathbb{Z}$. Then
\[t^m \zeta_{A[m]} = \zeta_A.\]
\item If $k'$ is the finite field extension of $k$ defined in (\ref{preg}) then $A' := A \otimes_k k'$ lies in $\grmod{k'H}$ and
\[\zeta_{A'} = \zeta_A.\]
\end{enumerate}
\end{lem}

\subsection{Rationality of graded Brauer characters}
\label{BrRat} By picking a homogeneous generating set we see that
any $M \in \mgr{R}$ actually lies in $\grmod{kH}$ and as such has a
graded Brauer character $\zeta_M$. It is easy to see that
$\zeta_M(1)$ is just the Poincar\'e series $P_M$ and is as such a
rational function (\ref{Dims}). The following result shows that
$\zeta_M(h)$ is also a rational function for any $h \in H_{\reg}$.
Let $m$ denote the $p'$-part of $|H|$.

\begin{thm} For any $M \in \mgr{R}$ and any $h \in H_{\reg}$ there exists a Laurent polynomial $u_h(t) \in F[t,t^{-1}]$ such that
\[\zeta_M(h) = \frac{u_h(t)}{(1 - t^m)^{d(M)}}.\]
\end{thm}
\begin{proof}
Without loss of generality we may assume that $H=\langle h\rangle$.
Moreover, since extension of scalars doesn't change the graded
Brauer character by Lemma \ref{GrBrProps}(iv), we may also assume
that $k = k'$. Hence $h$ acts diagonalizably on $V$.

We now prove the result by induction on $d(M) + \dim_k V$. Suppose first of all that $d(M)=0$. It follows that $M_n=0$ for sufficiently large $|n|$ and so $\zeta_M(h)$ is a Laurent polynomial as required. Since $d(M) \leq \dim_k V$ this also deals with the case when $\dim_k V = 0$, so we assume $\dim_k V > 0$ and $d(M) > 0$.

Choose an $h$-eigenvector $v \in V$ with eigenvalue $\lambda$. As $\lambda^m = 1$ we see that $h^{-1} v^m h = (v^h)^m = (\lambda v)^m = v^m$, so $z:= v^m$ is central in $R$.

Consider the graded submodule $T=\{\alpha\in M : \alpha.z^r=0$ for some $r\geq 0\}$ of $M$. Because $R$ is Noetherian and $M$ is finitely generated, $T$ is finitely generated as an $R$-module and as such is killed by some power of the central element $z$.

Choose an $h$-invariant complement $W$ for $kv$ in $V$ so that $\Sym(V) \cong \Sym(W)[v]$. It is now easy to see that $T$ is a finitely generated over $\Sym(W)$ and in fact $T \in \mgr{\Sym(W)\#H}$. Note that the dimension of $T$ viewed as a $\Sym(V)$-module is the same as the dimension of $T$ viewed as an $\Sym(W)$-module as both depend only on the Poincar\'e series of $T$.

Since $d(T)\leq d(M)$ and $\dim_k W < \dim_k V$, we know by induction that
\[\zeta_T(h)\cdot (1-t^m)^{d(M)} \in F[t,t^{-1}].\]
By Lemma \ref{GrBrProps}(i), $\zeta_M=\zeta_T+\zeta_{M/T}$, so it
now suffices to prove the result for the graded module $M/T$. So, by
replacing $M$ by $M/T$ we may assume $M$ is $z$-torsion-free.

Now as $z$ is a central element of $R$ of degree $m$, multiplication by $z$ induces a short exact sequence of graded $R$-modules
\[ 0\rightarrow M\stackrel{z}{\rightarrow} M[m]\rightarrow L[m]\rightarrow 0, \]
where $L:= M/Mz$. It follows from Lemma \ref{GrBrProps} that
\[t^m\zeta_M(h)=\zeta_M(h)+t^m\zeta_L(h).\]
But $L$ is a finitely generated graded $\Sym(W)\#H$-module and $d(L) \leq d(M) - 1$ by \cite[Proposition 11.3]{AM}, so by induction
\[\zeta_L(h)\cdot(1-t^m)^{d(M)-1}\in F[t,t^{-1}].\]
The result follows.
\end{proof}

Inspecting the proof shows that the Theorem is still valid with $m$
replaced by the order of $h$.

\section{Koszul resolutions}
We continue with the notation established in \S \ref{GrMod}. Let $d = \dim_k V$.

\subsection{The Koszul complex for graded $R$-modules}
\label{Kos}
With any finitely generated graded right $R$-module $M$ we associate the \emph{Koszul complex}
\[\mathbb{K}(M)_\ast := 0 \to M\otimes_k (\Lambda^dV)[-d] \stackrel{\phi_d}{\to} \cdots \stackrel{\phi_2}{\to} M \otimes_k (\Lambda^1 V)[-1] \stackrel{\phi_1}{\to} M \to 0\]
whose maps $\phi_j : M \otimes_k (\Lambda^jV)[-j] \to M\otimes_k (\Lambda^{j-1}V)[-j+1]$ are given by the usual formula
\[
\phi_j(m \otimes v_1 \wedge \cdots \wedge v_j) = \sum_{i=1}^j (-1)^{i+1} m.v_i \otimes  v_1 \wedge \cdots \wedge \hat{v}_i \wedge \cdots \wedge v_j,
\]
for any $m \in M$ and $v_1,\ldots, v_j \in V$. The square brackets indicate that we are thinking of $\Lambda^jV[-j]$ as a graded $k$-vector space concentrated in degree $j$.

\begin{lem} $\mathbb{K}(M)_\ast$ is a complex inside the abelian category $\grmod{kH}$.
\end{lem}
\begin{proof} It is well-known (and easily verified) that $\mathbb{K}(M)_\ast$ is a complex of $k$-vector spaces. Letting $H$ act diagonally on $M \otimes_k \Lambda^jV$ makes each $\phi_j$ into a map of right $kH$-modules as $h\cdot v_i^h = v_i\cdot h$ inside the ring $R$. Because of the shifts, each $\phi_j$ is also a map of graded modules of degree zero, as required.
\end{proof}

\subsection{Homology of $\mathbb{K}(M)_\ast$}
\label{KoszHom}
The projection map $\epsilon : R \twoheadrightarrow R_0$ with kernel $\bigoplus_{n=1}^\infty R_n$ is an algebra homomorphism from $R$ to $kH$ which gives $kH$ an $R$-$R$-bimodule structure. The Koszul complex is useful because it allows us to compute $\Tor_j^R(M,kH)$ for any $M\in \mgr{R}$.

\begin{prop}
\begin{enumerate}[{(}i{)}]
\item $\mathbb{K}(R)_\ast$ is a complex of $R$-$kH$-bimodules which is exact everywhere except in degree zero.
\item $H_0(\mathbb{K}(R)_\ast) \cong kH$ as $R$-$kH$-bimodules.
\item For all $M \in \mgr{R}$ and all $j\geq 0$ there is an isomorphism of $kH$-modules
\[H_j(\mathbb{K}(M)_\ast) \cong \Tor_j^R(M,kH).\]
\end{enumerate}
\end{prop}
\begin{proof}
(i) Consider first the special case when $H = 1$. Thinking of $V$ as an abelian $k$-Lie algebra, $R = \Sym(V)$ becomes the enveloping algebra of $V$ and
\[\mathbb{K}(R)_\ast = \Sym(V) \otimes_k \Lambda^\ast V\]
is just the Chevalley-Eilenberg complex of left $\Sym(V)$-modules \cite[\S 7.7]{Wei}. By \cite[Theorem 7.7.2]{Wei} $\mathbb{K}(R)_\ast$ is exact everywhere except in degree zero.

Returning to the general case, there is a natural isomorphism of complexes of $k$-vector spaces
\[R \otimes_{\Sym(V)} (\Sym(V)\otimes_k \Lambda^\ast V) \stackrel{\cong}{\longrightarrow} \mathbb{K}(R)_\ast.\]
Because $R$ is free of finite rank as a right $\Sym(V)$-module, it follows that $\mathbb{K}(R)_\ast$ is also exact everywhere except in degree zero. We use this isomorphism to give $\mathbb{K}(R)_\ast$ the structure of a complex of left $R$-modules.

On the other hand, $\mathbb{K}(R)_\ast$ is a complex of right $kH$-modules by Lemma \ref{Kos}. It can be checked that the two structures are compatible, so $\mathbb{K}(R)_\ast$ is a complex of $R$-$kH$-bimodules. Explicitly, the bimodule structure is given by the following formula:
\[s.(r\otimes v_1\wedge \cdots \wedge v_j).h = srh \otimes v_1^h \wedge \cdots \wedge v_j^h\]
for all $s,r \in R$, $v_1,\ldots,v_j \in V$ and $h \in H$.

(ii) The map $\epsilon : R \to kH$ which gives $kH$ its $R$-$kH$-bimodule structure is the cokernel of the first map $\phi_1 : R\otimes_k V \to R$ of the complex $\mathbb{K}(R)_\ast$. Hence
\[H_0(\mathbb{K}(R)_\ast) = R / \im \phi_1 \cong kH\]
as $R$-$kH$-bimodules.

(iii) Each term in $\mathbb{K}(R)_\ast$ is free of finite rank as a left $R$-module, so by parts (i) and (ii) $\mathbb{K}(R)_\ast$ is a free resolution of the left $R$-module $kH$. We can therefore use it to compute $\Tor_j^R(M,kH)$. Finally, the natural map $M \otimes_R \mathbb{K}(R)_\ast \to \mathbb{K}(M)_\ast$ is actually an isomorphism of complexes of right $kH$-modules, so
\[H_j(\mathbb{K}(M)_\ast) \cong H_j(M \otimes_R \mathbb{K}(R)_\ast) \cong \Tor_j^R(M,kH)  \]
as required.
\end{proof}

\subsection{A formula involving graded Brauer characters}
\label{KeyFormula} Let $M \in \mgr{R}$. By Lemma \ref{Kos} and
Proposition \ref{KoszHom}(iii), $\Tor_j^R(M,kH)$ is an object in
$\grmod{kH}$ and as such has a graded Brauer character
$\zeta_{\Tor_j^R(M,kH)}$. On the other hand, since $R$ is Noetherian
$M$ has a projective resolution consisting of finitely generated
$R$-modules. Computing $\Tor_j^R(M,kH)$ using this resolution shows
that these $kH$-modules are finite dimensional over $kH$, so
$\zeta_{\Tor_j^R(M,kH)}(h)$ is actually a Laurent polynomial in
$F[t,t^{-1}]$ for all $h\in H_{\reg}$.

\begin{prop} For any $h \in H_{\reg}$ and $M \in \mgr{R}$,
\[\sum_{j=0}^d (-1)^j \zeta_{\Tor_j^R(M,kH)}(h) = \zeta_M(h) \cdot \sum_{j=0}^d (-t)^j \varphi_{\Lambda^jV}(h).\]
\end{prop}
\begin{proof}In the notation of (\ref{HomEul}), Lemma \ref{GrBrProps}(i) says that $\psi_h : A \mapsto \zeta_A(h)$ is an additive function from the abelian category $\grmod{kH}$ to $F[[t,t^{-1}]]$ thought of as an abelian group. Applying Proposition \ref{KoszHom}(iii) and Lemma \ref{HomEul} we obtain
\[\sum_{j=0}^d (-1)^j \zeta_{\Tor_j^R(M,kH)}(h) = \psi_h(\Tor_\ast^R(M,kH)) = \psi_h(\mathbb{K}(M)_\ast).\]
Now $\mathbb{K}(M)_j = M \otimes_k \Lambda^jV[-j]$ so the result
follows from Lemma \ref{GrBrProps}(ii) and (iii).
\end{proof}

\section{Proof of Theorem \ref{MainResult}}
\label{PhiZeroProof}

\subsection{Another expression for $\Psi(g)$}\label {Schur}
Set $V := k\otimes_{\Fp} (N/N^p)$. Because $N$ is uniform, $N/N^p$ is an $\Fp$-vector space of
dimension $d$ and the right conjugation action of $G$ on $N$ induces
a right action of $G$ on $N/N^p$ by linear automorphisms. In this
way $V$ becomes a right $kG$-module with $\dim_kV = d$.

\begin{lem} For any $g \in G_{\reg}$, $\Psi(g)  =
\sum_{j=0}^d(-t)^j\varphi_{\Lambda^jV}(g).$
\end{lem}
\begin{proof}
Let $\beta = \Ad(g^{-1})\in\GL(\mathcal{L}(G))$ and let $j\geq 0$. Recall from (\ref{AdRep}) that $\Psi(g) = \det(1 - t\beta)$.
Now, $\Lambda^j N$ is a $\Lambda^j \beta$-stable lattice inside $\Lambda^j\mathcal{L}(G)$ whose reducion modulo $p$ is isomorphic to $\Lambda^j(N/N^p)$. Moreover, the endomorphism of $\Lambda^jV \cong \Lambda^j(N/N^p) \otimes_{\Fp} k$ induced by $\Lambda^j\beta$ is equal to the right action of $g$ on $\Lambda^jV$. The result now follows from the well-known formula
\[\det( 1 - t\beta ) = \sum_{j=0}^d (-t)^j\Tr(\Lambda^j \beta).\]
See, for example, \cite[p.487]{Serre3} or \cite[p. 77, (6.2)]{FulHar}.
\end{proof}

\begin{cor} The restriction of $\Psi$ to $G_{\reg}$ is locally constant.
\end{cor}
\begin{proof} Note that $V$ is a $k\overline{G}$-module because $[N,N]\leq N^p$ as $N$ is uniform, so the Brauer characters $\varphi_{\Lambda^jV}$ are constant on the cosets of $N$. Now apply the Lemma.
\end{proof}

\subsection{The associated graded ring} \label{AssocGr}
We now make the connection with the theory developed on the
preceeding two sections. Let $H := \overline{G} = G/N$; then $V$ is
a $kH$-module and we may form the skew group ring $\Sym(V) \# H$.
Recall from (\ref{GrBr}) that $w_N$ denotes the augmentation ideal
$(N-1)kG$.

\begin{lem}
The associated graded ring of $kG$ with respect to the $w_N$-adic
filtration is isomorphic to $R = \Sym(V) \# H$.
\end{lem}
\begin{proof}
When $N = G$ this is follows from \cite[Theorem 7.24]{DDMS}; see
also \cite[Lemma 3.11]{A}. Letting $\mathfrak{m}$ denote the
augmentation ideal of $kN$ we see that
\[\gr_{w_N} kG \cong kG \otimes_{kN}\gr_{\mathfrak{m}}kN \cong kH\otimes_k\Sym(V)\]
as a right $\Sym(V)$-module, because $kN$ acts trivially on its graded ring $\Sym(V)$. Moreover, the zero$^{th}$ graded part of $\gr_{w_N}kG$ is isomorphic to $kH$ as a $k$-algebra, so $H$ embeds into the group of units of $\gr_{w_N} kG$. It is now easy to verify that the multiplication works as needed.
\end{proof}

We will identify $R$ with $\gr_{w_N}kG$ in what follows.

\subsection{A spectral sequence}
\label{SpecSeq} The last step in the proof involves relating $\Tor$
groups over $kG$ with $\Tor$ groups over the associated graded ring
$R$. There is a standard spectral sequence originally due to Serre
which does the job.

\begin{prop} For any finitely generated $kG$-module $M$ there exists a homological spectral sequence in $\mathcal{M}(kH)$
\[E^1_{ij} = \Tor_{i+j}^R(\gr M, kH)_{\degree -i} \Longrightarrow \Tor_{i+j}^{kG}(M,kH).\]
\end{prop}
\begin{proof} As in (\ref{GrBr}), we only consider the \emph{deduced filtration} on $M$, given by $M^n = Mw_N^n$ for $n\geq 0$. As $M$ is finitely generated over $kG$, the associated graded module
$\gr M$ is finitely generated over the Noetherian ring $\gr kG\cong
R$. By the proof of \cite[Proposition 3.4]{AB}, the $w_N$-adic
filtration on $kG$ is complete.

We claim that $M$ is complete with respect to the deduced
filtration. Because $kG$ is Noetherian we can find an exact sequence
$(kG)^a \stackrel{\alpha}{\to} (kG)^b \stackrel{\beta}{\to} M \to 0$
in $\mathcal{M}(kG)$. Giving all the modules involved deduced
filtrations, $(kG)^a$ and $(kG)^b$ are compact and the maps
$\alpha,\beta$ are continuous. Hence $\im \alpha$ and $M$ are
compact, so $\im \alpha$ is closed in $(kG)^b$ and $M$ is complete
as claimed.

We can now apply \cite[Proposition 8.1]{Totaro} to the modules $A =
M$ and $B = kH$ over the complete filtered $k$-algebra $kG$, where
we equip $B$ with the trivial filtration $B^0 = B$ and $B^1=0$. This
gives us the required spectral sequence of $k$-vector spaces.
Examining the construction shows that it is actually a spectral
sequence in $\mathcal{M}(kH)$.
\end{proof}

\subsection{Proof of Theorem \ref{MainResult}}
\label{ProofMainRes}
The spectral sequence $E$ of Proposition \ref{SpecSeq} gives us suitable filtrations on the $kH$-modules $\Tor_n^{kG}(M,kH)$. In the notation of (\ref{SpecEul}) we can rewrite the information we gain from the spectral sequence as follows:
\[\Tot(E^1)_n = \Tor_n^R(\gr M, kH) \quad\mbox{and} \quad \Tot(E^\infty)_n = \gr \Tor_n^{kG}(M,kH)\]
for each $n \in \mathbb{Z}$. We have already observed in
(\ref{KeyFormula}) that $\Tor_n^R(\gr M, kH)$ is a finite
dimensional $kH$-module, which is moreover zero whenever $n > d$ by
Proposition \ref{KoszHom}(iii). Thus $E$ is totally bounded. Now,
Lemma \ref{BrCh}(ii) shows that $A \mapsto \varphi_A(g)$ is an
additive $F$-valued function on the objects of $\mathcal{M}(kH)$. By
Proposition \ref{SpecEul},
\[\rho_N[M](g) = \sum_{j=0}^d (-1)^j \varphi_{\Tor_j^{kG}(M,kH)}(g) = \sum_{j=0}^d (-1)^j \varphi_{\Tor_j^R(\gr M,kH)}(g).\]

We may now apply Proposition \ref{KeyFormula} and Lemma \ref{Schur}
to obtain
\[\sum_{j=0}^d (-1)^j \varphi_{\Tor_j^R(\gr M,kH)}(g) = \sum_{j=0}^d (-1)^j \zeta_{\Tor_j^R(\gr M,kH)}(g) |_{t=1} =
\left(\zeta_{\gr M}(g) \cdot \Psi(g)\right) |_{t=1},\] as
required.\qed

\section{Euler characteristics}
\subsection{Twisted $\mu$-invariants}
\label{EulerCharMu} Because we only deal with Iwasawa modules which
are killed by $p$ in this paper, it is easy to see that the
definition of the \emph{Euler characteristic} \cite[\S 1.5]{AW} of a
finitely generated $kG$-module $M$ of finite projective dimension
can be given as follows:
\[\chi(G,M) := \prod_{j\geq 0} |\Tor^{kG}_j(M,k)|^{(-1)^j}.\]

Let $\{V_1,\ldots,V_s\}$ be a complete list of representatives for
the isomorphism classes of simple $kG$-modules. The \emph{$i$-th
twisted $\mu$-invariant} of $M$ \cite[\S 1.5]{AW} for $i=1,\ldots,
s$ is defined by the formula
\[\mu_i(M) = \frac{\log_q \chi(G, M \otimes_k V_i^*)}{\dim_k \End_{kG}(V_i)}, \]
where $V_i^*$ is the dual module to $V_i$. We assume that $V_1$ is
the trivial $kG$-module $k$, so that
\[\mu_1(M) = \log_q \chi(G,M).\]

We proved in \cite{AW} that these twisted $\mu$-invariants
completely determine the characteristic element of $M$ viewed as an
$\mathcal{O}G$-module \cite[Theorem 1.5]{AW}. This adds to the motivation of the problem of computing the Euler characteristic $\chi(G,M)$.

\subsection{The base change map}
\label{kGbar} Before we can proceed, we need to record some
information about the base change map $\theta_N : \mathcal{G}_0(kG)
\to \mathcal{G}_0(k\overline{G})$.

We say that a map $f$ between two abelian groups is an
\emph{$\mathbb{Q}$-isomorphism}
 if it becomes an isomorphism after tensoring with $\mathbb{Q}$. Equivalently, $f$ has torsion kernel and cokernel.

\begin{prop} There is a commutative diagram of Grothendieck groups
\[
\xymatrix{
K_0(kG) \ar[d]_{c}\ar[r]^{\pi_N} & K_0(k\overline{G}) \ar[d]^{c_N} \\
\mathcal{G}_0(kG) \ar[r]_{\theta_N} & \mathcal{G}_0(k\overline{G}). \\
    }
\]
The map $\pi_N$ is an isomorphism and the other maps are $\mathbb{Q}$-isomorphisms.
\end{prop}
\begin{proof}
The vertical maps $c$ and $c_N$ are called \emph{Cartan maps} and
are defined by inclusions between admissible subcategories. The base
change map $\pi_N$ is defined by $\pi_N[P] =
[P\otimes_{kG}k\overline{G}]$ for all $P \in \mathcal{P}(kG)$, and
we have already discussed $\theta_N$ in (\ref{BaseChange}). This
well-known diagram appears in \cite[p. 454]{Bass} and expresses the
fact that the Cartan maps form a natural transformation from
$K$-theory to $\mathcal{G}$-theory.

Now, $kG$ is a \emph{crossed product} of $kN$ with the finite group
$\overline{G}$:
\[kG \cong kN \ast \overline{G},\]
see, for example \cite[\S 2.3]{AB2} for more details. Because $kN$
is a Noetherian $k$-algebra of finite global dimension, the map $c$
is an $\mathbb{Q}$-isomorphism by a general result on the $K$-theory of
crossed products \cite{AW3}. Considering the case when $G$ is finite shows that $c_N$ is an $\mathbb{Q}$-isomorphism as well -- this also follows from a well-known result of Brauer: see, for example, \cite[Corollary 1 to Theorem 35]{Serre2}.

Finally $\pi_N$ is an isomorphism because $kG$ is a complete semilocal
ring -- see \cite[Lemma 2.6 and Proposition 3.3(a)]{AW}. It follows
that $\theta_N$ must also be an $\mathbb{Q}$-isomorphism, as required.
\end{proof}

Recalling Theorem \ref{BermanWitt} and Lemma \ref{Dev}, we obtain

\begin{cor} The map $\rho_N$ featuring in Theorem \ref{MainResult} is a $\mathbb{Q}$-isomorphism.
\end{cor}

If the group $G$ has no elements of order $p$, then $kG$ has finite
global dimension and the Cartan map $c$ is actually an isomorphism
by Quillen's Resolution Theorem \cite[Theorem 12.4.8]{MCR}. In this
case, therefore, we do not have to rely on \cite{AW3}.

\subsection{Computing Euler characteristics using $\theta_N$}
\label{pairing}
Let $P \in \mathcal{P}(k\overline{G})$ and let $V \in \mathcal{M}(k\overline{G})$. The rule
\[(P,V) \mapsto \dim_k \Hom_{k\overline{G}}(P,V)\]
defines an additive function from $\mathcal{P}(k\overline{G})\times \mathcal{M}(k\overline{G})$ to $\mathbb{Z}$ and hence a pairing
\[\langle-,-\rangle_N:K_0(k\overline{G})\times\mathcal{G}_0(k\overline{G})\rightarrow\mathbb{Z}.\]
This pairing appears in \cite[p. 121]{Serre2}. Now, by extending scalars, we can define a bilinear form
\[\langle-,-\rangle_N : \mathbb{Q}K_0(k\overline{G})\times \mathbb{Q}\mathcal{G}_0(k\overline{G}) \to \mathbb{Q}\]
which is in fact non-degenerate. We saw in Proposition \ref{kGbar} that the Cartan map $c_N : \mathbb{Q}K_0(k\overline{G}) \to \mathbb{Q}\mathcal{G}_0(k\overline{G})$ is an isomorphism. This allows us to define a non-degenerate bilinear form
\[(-,-)_N : \mathbb{Q}\mathcal{G}_0(k\overline{G})\times \mathbb{Q}\mathcal{G}_0(k\overline{G}) \to \mathbb{Q}\]
by setting $(x,y)_N = \langle c_N^{-1}(x),y\rangle_N$ for $x,y \in \mathbb{Q}\mathcal{G}_0(k\overline{G})$.

\begin{prop} For any finitely generated $kG$-module $M$ of finite projective dimension, the Euler characteristic of $M$ can be computed as follows:
\[\log_q\chi(G,M) = (\theta_N[M],k)_N.\]
\end{prop}
\begin{proof}
Suppose first that $M$ is projective. The usual adjunction between $\otimes$ and $\Hom$ gives isomorphisms
\[\Hom_k(M\otimes_{kG}k,k) \cong  \Hom_{kG}(M,k) \cong \Hom_{k\overline{G}}(M\otimes_{kG}k\overline{G},k).\]
As these $\Hom$ spaces are finite dimensional over $k$, we obtain
\[\dim_k (M\otimes_{kG}k)=\langle [M\otimes_{kG}k\overline{G}],k\rangle_N .\]
Now because $M$ is projective, $\Tor_i^{kG}(M,k) = 0 =
\Tor_i^{kG}(M,k\overline{G})$ for $i>0$, and $\theta_N[M] =
c_N(\pi_N[M])$ by Proposition \ref{kGbar}. Hence
\[\log_q\chi(G,M) = \dim_k (M\otimes_{kG}k) = \langle \pi_N[M],k \rangle_N = (\theta_N[M],k)_N\]
as required. Returning to the general case, if $0 \to P_n \to \cdots \to P_0 \to M \to 0$ is a resolution of $M$ in $\mathcal{P}(kG)$ then Lemma \ref{HomEul} gives
\[\theta_N[M] = \sum_{i=0}^n (-1)^i \theta_N[P_i]\quad\mbox{and}\quad \log_q\chi(G,M) = \sum_{i=0}^n(-1)^i\log_q\chi(G,P_i).\]
The result now follows from the first part.
\end{proof}

\subsection{Euler characteristics for modules of infinite projective dimension}
\label{GenEulCh} The definition of $\chi(G,M)$ given in
(\ref{EulerCharMu}) only makes sense when the module $M$ has finite
projective dimension. However, the expression $(\theta_N[M],k)_N$ makes sense for arbitrary $M \in \mathcal{M}(kG)$. We include the subscripts, because \emph{a priori} this depends on the choice of the
open normal uniform subgroup $N$.

\begin{lem} Let $M$ be a finitely generated $kG$-module. Then
\begin{enumerate}[{(}i{)}]
\item $\psi_N : M \mapsto (\theta_N[M],k)_N$ is an additive function on the objects of $\mathcal{M}(kG)$,
\item $\psi_N(M)$ does not depend on the choice of $N$.
\end{enumerate}
\end{lem}
\begin{proof}
It suffices to prove part (ii). Now if $U\in\mathbf{U}_{G,p}$ is another uniform subgroup of $G$, then $\psi_U(M) = \log_q \chi(G,M) = \psi_N(M)$ whenever $M$ is projective, by Proposition
\ref{EulerCharMu}. Rephrasing this in the language of Grothendieck
groups,
\[\mathbb{Q}K_0(kG) \stackrel{c}{\longrightarrow} \mathbb{Q}\mathcal{G}_0(kG) \stackrel{\psi_U - \psi_N}{\longrightarrow} \mathbb{Q}\]
is a complex of $\mathbb{Q}$-vector spaces. By Proposition \ref{kGbar}, the
Cartan map $c$ is an isomorphism, so $\psi_U(M) = \psi_N(M)$ for \emph{any}
finitely generated $kG$-module $M$, as required.
\end{proof}

In view of this result, we propose to extend the definition given in
(\ref{EulerCharMu}) as follows.

\begin{defn} The \emph{Euler characteristic} of a finitely generated $kG$-module $M$ is defined to be
\[\chi(G,M) := q^{(\theta_N[M],k)_N} \in q^{\mathbb{Q}}\]
for any choice of open normal uniform subgroup $N$ of $G$.
\end{defn}
\subsection{Trivial Euler characteristics}
\label{TrivEulCh} First, a preliminary

\begin{lem} For any $g\in G_{\reg}$, the multiplicity of $1$ as a root of the polynomial $\Psi(g)$ equals $\dim C_G(g)$.
\end{lem}
\begin{proof} As $\Psi(g)\cdot\det\Ad(g)$ is the characteristic polynomial of $\Ad(g)$, the first number equals the dimension of the space $C:=\{x \in \mathcal{L}(G) : \Ad(g)(x) = x\}$. The definition of $\Ad(g)$ shows that $C\cap N$ is just the centralizer $C_N(g)$ of $g$ in $N$. Because $C_N(g) = C_G(g)\cap N$ is open in $C_G(g)$, the dimension of $C_G(g)$ as a compact $p$-adic analytic group equals the dimension of $C$ as a $\Qp$-vector space. The result follows.
\end{proof}

Recall \cite[\S 5.4]{AB2} that as $kG$ is an Auslander-Gorenstein ring, every finitely generated $kG$-module $M$ has a \emph{canonical dimension} which we will denote by $d(M)$. By \cite[\S 5.4(3)]{AB2} this is a non-negative integer which equals the dimension $d(\gr M)$ of the associated graded module $\gr M$, defined in (\ref{Dims}).

\begin{prop} Let $M$ be a finitely generated $kG$-module and let $g \in G_{\reg}$ be such that $d(M) < \dim C_G(g)$. Then $\rho_N[M](g)=0$.
\end{prop}
\begin{proof}
By Theorem \ref{MainResult} and Theorem \ref{BrRat} there exists a Laurent polynomial $u_g(t) \in F[t,t^{-1}]$ such that $\rho_N[M](g)$ equals the value at $t=1$ of the rational function
\[ \frac{ u_g(t) \cdot \Psi(g) } { (1 - t^m)^{d(\gr M)} }.\]
Because we are assuming that $\dim C_G(g) > d(M) = d(\gr M)$, this rational function has a zero at $t=1$ in view of the Lemma.
\end{proof}

We can now give our first application of Theorem \ref{MainResult}.

\begin{proof}[Proof of Theorem \ref{TrivialEul}] By the Proposition, $\rho_N[M]=(\varphi \circ \lambda_N)(\theta_N[M])= 0$. As $\varphi$ and $\lambda_N$ are isomorphisms by Theorem \ref{BermanWitt} and Lemma \ref{Dev}, $\theta_N[M] = 0$. The result now follows from the new definition of $\chi(G,M)$ given in (\ref{GenEulCh}).
\end{proof}

\section{$K$-theory}
\label{Kth}
\subsection{The dimension filtration}\label{DimFilt}
Recall from (\ref{BrCh}) that $\mathcal{F}_0$ denotes the category of all $kG$-modules which are finite dimensional over $k$.

Now, a finitely generated $kG$-module $M$ is finite dimensional over $k$ if and only if $d(M) = d(\gr M) = 0$, because both conditions are equivalent to the Poincar\'e series (\ref{Dims}) of $\gr M$ being a polynomial in $t$. We can therefore unambiguously define $\mathcal{F}_i = \mathcal{F}_i(G)$ to be full subcategory of $\mathcal{M}(kG)$ consisting of all modules $M$ with $d(M) \leq i$, for each $i=0,\ldots,d$. Thus we have an ascending chain of subcategories
\[ \mathcal{F}_0 \subseteq \mathcal{F}_1 \subseteq \cdots \subseteq \mathcal{F}_{d-1} \subseteq \mathcal{F}_d = \mathcal{M}(kG).\]

Using \cite[\S 5.3]{AB2} we see that each $\mathcal{F}_i$ is an admissible (in fact, Serre) subcategory of $\mathcal{M}(kG)$, so we can form the Grothendieck groups $K_0(\mathcal{F}_i)$. The inclusions $\mathcal{F}_i \subseteq \mathcal{F}_d$ induce maps
\[\alpha_i : K_0(\mathcal{F}_i) \to K_0(\mathcal{F}_d) \qquad \mbox{and}\qquad \alpha_i : FK_0(\mathcal{F}_i) \to FK_0(\mathcal{F}_d) \]
which we would like to understand. Because $\rho_N$ is a
$\mathbb{Q}$-isomorphism by Corollary \ref{kGbar}, we focus on the
image of $\rho_N \circ \alpha_i$.

\subsection{Module structures}
\label{ModStr}
Note that $K_0(\mathcal{F}_0)$ is a commutative ring with multiplication induced by the tensor product. By adapting the argument used in \cite[Proposition 7.3]{AW} and using \cite[\S 5.4(5)]{AB2}, we see that the twist $M\otimes_k V$ of any $M \in \mathcal{F}_d$ and any $V \in \mathcal{F}_0$ satisfies $d(M \otimes_k V) = d(M)$. In this way $K_0(\mathcal{F}_i)$ becomes a $K_0(\mathcal{F}_0)$-module and it is clear that the maps $\alpha_i : K_0(\mathcal{F}_i) \to K_0(\mathcal{F}_d)$ respect this module structure.

\begin{lem} The map $\lambda_N \circ \theta_N : K_0(\mathcal{F}_d) \to K_0(\mathcal{F}_0)$ is a map of $K_0(\mathcal{F}_0)$-modules.
\end{lem}
\begin{proof} Let $V \in \mathcal{M}(k\overline{G})$. For any $M\in\mathcal{F}_d$ we can define a function
\[\beta_M : (M\otimes_k V) \otimes_{kG}k\overline{G} \to (M\otimes_{kG}k\overline{G})\otimes_k V\]
by the formula $\beta_M((m\otimes v)\otimes h) = (m\otimes h)\otimes vh$ for $m \in M$, $v \in V$ and $h \in \overline{G}$. A straightforward check shows that $\beta_M$ is a homomorphism of right $kG$-modules with inverse $\gamma_M$, defined by the formula $\gamma_M((m\otimes h)\otimes v) = (m\otimes vh^{-1}) \otimes h$. Hence the functors $M \mapsto (M\otimes_k V) \otimes_{kG}k\overline{G}$ and $M \mapsto (M\otimes_{kG}k\overline{G})\otimes_k V$ are isomorphic, so
\[\Tor_j^{kG}(M\otimes_k V,k\overline{G}) \cong \Tor_j^{kG}(M,k\overline{G})\otimes_k V\]
for all $j\geq 0$. It follows that $\theta_N[M\otimes_k V] = \theta_N[M].[V]$ for all $M\in \mathcal{F}_d$ and $V \in \mathcal{M}(k\overline{G})$. Because $\lambda_N : \mathcal{G}_0(k\overline{G}) \to K_0(\mathcal{F}_0)$ is an isomorphism by Lemma \ref{Dev},
\[\lambda_N(\theta_N[M\otimes_k V]) = \lambda_N(\theta_N[M]).[V]\]
for all $M\in \mathcal{F}_d$ and $V \in \mathcal{F}_0$, as required.
\end{proof}

Hence the image of $\lambda_N\circ\theta_N\circ\alpha_i$ is always an \emph{ideal} of $K_0(\mathcal{F}_0)$. On the other hand, $C(G_{\reg};F)^{G\times \mathcal{G}_k}$ is a commutative $F$-algebra via pointwise multiplication of functions, and Lemma \ref{BrCh}(iii) shows that the map
\[\varphi : FK_0(\mathcal{F}_0) \to C(G_{\reg};F)^{G\times \mathcal{G}_k}\]
appearing in Theorem \ref{BermanWitt} is an $F$-algebra isomorphism, so $\im (\rho_N \circ \alpha_i)$ is always an ideal of $C(G_{\reg};F)^{G\times \mathcal{G}_k}$.

It is easy to see that the ideals of this algebra are in bijection with the subsets of the orbit space $(G\times\mathcal{G}_k)\backslash G_{\reg}$; which subset does $\im (\rho_N \circ \alpha_i)$ correspond to?

\subsection{An upper bound for $\rk \alpha_i$}
\label{UpBndPf}
Define a subset $S_i$ of $G_{\reg}$ by the formula
\[S_i := \{g \in G_{\reg} : \dim C_G(g) \leq i\}.\]
We record some basic facts about these subsets of $G_{\reg}$.

\begin{lem}
 \begin{enumerate}[{(}i{)}]
\item $S_i$ is a union of conjugacy classes in $G$.
\item $S_i$ is stable under the action of $\mathcal{G}_k$ on $G_{\reg}$.
\item $S_i$ is a clopen subset of $G_{\reg}$ and hence a closed subset of $G$.
\end{enumerate}
\end{lem}
\begin{proof} (i) This is clear.

(ii) If $g$ is a power of $h$ then $C_G(h) \leq C_G(g)$, therefore
if $g$ and $h$ lie in the same $\mathcal{G}_k$-orbit inside
$G_{\reg}$ then their centralizers are equal.

(iii) This follows from Corollary \ref{HallBij}.
\end{proof}

We can now give our second application of Theorem \ref{MainResult}.
\begin{proof}[Proof of Theorem \ref{UpBnd}]
By Proposition \ref{TrivEulCh}, $\rho_N[M]$ is zero on $G_{\reg} - S_i$ for any $M \in \mathcal{F}_i$. Hence $\im (\rho_N \circ \alpha_i) \subseteq C(S_i;F)^{G\times\mathcal{G}_k}$ and the result follows from Corollary \ref{kGbar}.
\end{proof}

\section{Some special cases}

\subsection{The rank of $\alpha_d$}
\label{KoFd}
We begin by recording the rank of $\alpha_d$, or equivalently the rank of $K_0(\mathcal{F}_d)=\mathcal{G}_0(kG)$.
\begin{prop} The rank of $\alpha_d$ equals the number of $G\times \mathcal{G}_k$-orbits on $G_{\reg}$:
\[\rk \alpha_d = \rk K_0(\mathcal{F}_d) = |(G\times\mathcal{G}_k)\backslash G_{\reg}|.\]
\end{prop}
\begin{proof} This follows from Corollary \ref{kGbar}.
\end{proof}

\subsection{A localisation sequence}
\label{LocSeq}
Consider the localisation sequence of $K$-theory \cite[Theorem 5.5]{Q} for the Serre subcategory $\mathcal{F}_i$ of the abelian category $\mathcal{F}_d$:
\[K_0(\mathcal{F}_i) \stackrel{\alpha_i}{\longrightarrow} K_0(\mathcal{F}_d) \longrightarrow K_0(\mathcal{F}_d/\mathcal{F}_i)\longrightarrow 0.\]
Because we already know the rank of $K_0(\mathcal{F}_d)$, this
sequence shows the problem of computing the rank of $\alpha_i$ is
equivalent to the problem of computing the rank of the Grothendieck
group of the quotient category $\mathcal{F}_d/\mathcal{F}_i$. At
present, the only non-trivial case when we understand $\mathcal{F}_d
/\mathcal{F}_i$ is the case $i=d-1$ --- see (\ref{torsion}). First,
we require some results from noncommutative algebra.

\subsection{Artinian rings and minimal primes} \label{ArtMin}
Recall that if $R$ is a (not necessarily commutative) ring, then an ideal $I$ of $R$ is said to be \emph{prime} if whenever $A,B$ are ideals of $R$ strictly containing $I$, their product $AB$ also strictly contains $I$. A \emph{minimal prime} of $R$ is a prime ideal which is minimal with respect to inclusion - equivalently, it has height zero. The following result is well-known.

\begin{lem} Let $R$ be an Artinian ring. Then
\begin{enumerate}[{(}i{)}]
\item $\mathcal{G}_0(R)$ is a free abelian group on the isomorphism classes of simple $R$-modules.
\item There is a natural bijection between the isomorphism classes of simple $R$-modules and the minimal primes of $R$, given by
\[[M] \mapsto \Ann_R(M).\]
\end{enumerate}
\end{lem}

\subsection{The finite radical}
Recall \cite[1.3]{AB} the important characteristic subgroup $\Delta^+$ of $G$, defined by
\[\Delta^+ = \{x \in G : |G:C_G(x)| < \infty \quad \mbox{and} \quad o(x) < \infty\}.\]
This group is sometimes called the \emph{finite radical} of $G$ and consists of all elements of finite order in $G$ whose conjugacy class is finite. It is also the largest finite normal subgroup of $G$. In our notation, $\Delta^+_{\reg} = \Delta^+ \cap G_{\reg}$ is just the complement of $S_{d-1}$ in $G_{\reg}$:
\[\Delta^+_{\reg} = G_{\reg} - S_{d-1}\]
and as such is a union of $G\times\mathcal{G}_k$-orbits in $G_{\reg}$.

\subsection{The classical ring of quotients $Q(kG)$}
\label{QkG} By \cite[Proposition 7.2]{AW} $kG$ has an Artinian
classical ring of quotients $Q(kG)$. We have already computed the
rank of $\mathcal{G}_0(Q(kG))$ under the assumption that the order
of the finite group $\Delta^+$ is coprime to $p$ \cite[Theorem
12.7(b)]{AW}. We can now present a generalization of this result,
valid without any restrictions on $G$.

\begin{thm} The rank of $\mathcal{G}_0(Q(kG))$ equals the number of $G\times \mathcal{G}_k$-orbits on $\Delta^+_{\reg}$:
\[\rk \mathcal{G}_0(Q(kG)) = |(G\times\mathcal{G}_k)\backslash \Delta^+_{\reg}|.\]
\end{thm}
\begin{proof}
We will show that both numbers in question are equal to the number
of minimal primes of $kG$, $r$ say. Let $J$ denote the \emph{prime
radical} of $kG$, defined as the intersection of all prime ideals of
$kG$. By passing to $kG/J$ and applying \cite[Proposition
3.2.2(i)]{MCR}, we see that minimal primes of $kG$ are in bijection
with the minimal primes of $Q(kG)$. As $Q(kG)$ is Artinian, Lemma
\ref{ArtMin} implies that $\rk \mathcal{G}_0(Q(kG)) = r$.

Next, the group $G$ acts on $\Delta^+$ by conjugation and therefore
permutes the the minimal primes of $k\Delta^+$. It was shown in
\cite[Theorem 5.7]{A2} that there is a natural bijection between the
minimal primes of $kG$ and $G$-orbits on the minimal primes of
$k\Delta^+$. The group $G$ also permutes the simple
$k\Delta^+$-modules and respects the correspondence between these
and the minimal primes of $k\Delta^+$ given in Lemma
\ref{ArtMin}(ii). Hence $r$ is also the number of $G$-orbits on the
simple $k\Delta^+$-modules. Finally, the $G$-equivariant version of
the Berman--Witt Theorem \cite[Corollary 12.6]{AW} shows that the
latter is just $|(G\times\mathcal{G}_k)\backslash \Delta^+_{\reg}|$
and the result follows.
\end{proof}

\subsection{The rank of $\alpha_{d-1}$}
\label{torsion} Recall that a finitely generated $kG$-module is
\emph{torsion} if and only if $M \otimes_{kG} Q(kG) = 0$. By
\cite[Lemma 1.4]{CSS}, $M$ is torsion if and only if $d(M) < d(kG) =
d$, so $\mathcal{F}_{d-1}$ is just the category of all finitely
generated torsion $kG$-modules, as mentioned in the introduction.

\begin{lem}
The quotient category $\mathcal{F}_d / \mathcal{F}_{d-1}$ is
equivalent to $\mathcal{M}(Q(kG))$.
\end{lem}
\begin{proof} This follows from \cite[Propositions XI.3.4(a) and
XI.6.4]{Sten}, with appropriate modifications to handle the finitely
generated case.
\end{proof}

We can now use the localisation sequence of (\ref{LocSeq}) to show
that the upper bound for $\rk \alpha_i$ given in Theorem \ref{UpBnd}
is attained in the case when $i=d-1$.

\begin{prop} The rank of $\alpha_{d-1}$ equals the number of $G\times\mathcal{G}_k$-orbits on $S_{d-1}$:
\[\rk \alpha_{d-1} = |(G \times \mathcal{G}_k) \backslash S_{d-1}|.\]
\end{prop}
\begin{proof}
In view of the Lemma, the localisation sequence becomes
\[K_0(\mathcal{F}_{d-1}) \stackrel{\alpha_{d-1}}{\longrightarrow} K_0(\mathcal{F}_d) \longrightarrow \mathcal{G}_0(Q(kG)) \longrightarrow 0.\]
Hence $\rk (\alpha_{d-1})= \rk K_0(\mathcal{F}_d) - \rk
\mathcal{G}_0(Q(kG))$. Now apply Proposition \ref{KoFd} and Theorem
\ref{QkG}, bearing in mind that $S_{d-1} = G_{\reg} -
\Delta^+_{\reg}$.
\end{proof}

\subsection{The rank of $\alpha_0$}
\label{rkalpha0}
We will see in (\ref{CEx}) that the rank of $\alpha_i$ does \emph{not} always attain the upper bound of Theorem \ref{UpBnd}. Here is another special case when $\rk \alpha_i$ is well-behaved.

\begin{prop} The rank of $\alpha_0$ equals the number of $G\times\mathcal{G}_k$-orbits on $S_0$:
\[\rk \alpha_0 = |(G \times \mathcal{G}_k) \backslash S_0|.\]
\end{prop}
\begin{proof}
Let $M\in\mathcal{F}_0$. Because $M$ is finite dimensional over $k$, the graded Brauer character
$\zeta_{\gr M}$ is a polynomial, so
\[\varphi_M = \zeta_{\gr M}|_{t=1},\]
thought of as $F$-valued functions on $G_{\reg}$. Applying the
explicit formula for $\rho_N$ given in Theorem \ref{MainResult} shows
that
\[\rho_N[M](g) = \varphi_M(g) \cdot \det(1 - \Ad(g^{-1}))\]
for any $g \in G_{\reg}$. We therefore have a commutative diagram
\[
\xymatrix{
FK_0(\mathcal{F}_0) \ar[r]^{\varphi}\ar[d]_{\alpha_0} & C(G_{\reg};F)^{ G\times\mathcal{G}_k} \ar[d]^{\eta}\\
FK_0(\mathcal{F}_d) \ar[r]_{\rho_N} & C(G_{\reg};F)^{ G\times\mathcal{G}_k}
}
\]
where $\eta$ is multiplication by the locally constant function
\[\Psi|_{t=1} : g \mapsto \det(1 - \Ad(g^{-1})).\]
Using Lemma \ref{TrivEulCh} we see that $\Psi|_{t=1}(g) \neq 0$ if and only if $\dim C_G(g) = 0$. It follows that the image of $\eta$ is precisely $C(S_0;F)^{G\times\mathcal{G}_k}$. As the maps $\varphi$ and $\rho_N$ are isomorphisms by Theorem \ref{BermanWitt} and Corollary \ref{kGbar},
\[\rk \alpha_0 = \rk \eta = |(G\times\mathcal{G}_k)\backslash S_0|\]
as required.
\end{proof}

We now start preparing for the proof of Theorem \ref{VirtAb} which says that the upper bound of Theorem \ref{UpBnd} is always attained if the group $G$ is virtually abelian.

\section{Induction of modules}

\subsection{Dimensions}
\label{IndMod}
In what follows we fix a closed subgroup $H$ of $G$ of dimension $e$. Recall \cite[Lemma 4.5]{Bru} that $kG$ is a \emph{flat} $kH$-module. Therefore the induction functor
\[\Ind_H^G : \mathcal{M}(kH) \to \mathcal{M}(kG)\]
which sends $M \in \mathcal{M}(kH)$ to $M\otimes_{kH}kG\in\mathcal{M}(kG)$ is exact and induces a map
\[\Ind_H^G : \mathcal{G}_0(kH) \to \mathcal{G}_0(kG).\]
We can obtain precise information about the dimension of an induced module.
\begin{lem} Let $M\in\mathcal{M}(kH)$. Then
\[d(M\otimes_{kH}kG) = d(M) + d - e.\]
\end{lem}
\begin{proof}
Recall that the \emph{grade} $j(X)$ of a finitely generated $R$-module $X$ over an Auslander-Gorenstein ring $R$ is defined by the formula
\[j(X) = \min\{j : \Ext^j_R(X, R) \neq 0\}.\]
The \emph{canonical dimension} of $X$ is the non-negative integer $\id(R) - j(X)$ where $\id(R)$ is the injective dimension of $R$.

Now, choosing a free resolution of $M$ and using the fact that $kG$ is a flat $kH$-module, we see that
\[kG\otimes_{kH} \Ext_{kH}^j(M,kH) \cong \Ext_{kG}^j(M\otimes_{kH}kG,kG)\]
as left $kG$-modules, for any $j\geq 0$. In fact, $kG$ is a \emph{faithfully flat} (left) $kH$-module by \cite[Lemma 5.1]{A2}, which implies that $kG \otimes_{kH}A = 0$ if and only if $A = 0$ for any finitely generated left $kH$-module $A$.

Hence $j(M) = j(M\otimes_{kH}kG)$ and the result follows because $\id(kG) = \dim G = d$ and $\id(kH) = \dim H = e$.
\end{proof}

\subsection{Proposition}\label{IndComDiag} Suppose that the group $H\cap N$ is uniform. Then there exists a map $\iota_N : C(H_{\reg};F)^{H\times\mathcal{G}_k} \to C(G_{\reg};F)^{G\times\mathcal{G}_k}$ such that the following diagram commutes:
\[
\xymatrix{
FK_0(\mathcal{F}_{d-e}(G)) \ar[r]^{\alpha_{d-e}} & FK_0(\mathcal{F}_d(G)) \ar[r]^{\rho_N} & C(G_{\reg};F)^{G\times\mathcal{G}_k}\\
FK_0(\mathcal{F}_0(H)) \ar[u]^{\Ind_H^G}\ar[r]_{\alpha_0} & FK_0(\mathcal{F}_e(H)) \ar[u]^{\Ind_H^G}\ar[r]_{\rho_{H\cap N}} & C(H_{\reg};F)^{H \times \mathcal{G}_k} \ar[u]_{\iota_N}.
}
\]
\begin{proof} We assume that $H\cap N$ is uniform only to make sure that the map $\rho_{H\cap N}$ makes sense. We construct this diagram in several steps --- note that the left-hand square makes sense by Lemma \ref{IndMod}.

Let $\overline{H} := HN/N \cong H/(H\cap N)$, let $\overline {G} := G/N$ and define a map
\[\Ind_{\overline{H}}^{\overline{G}} : C(\overline{H}_{\reg};F)^{\overline{H}\times\mathcal{G}_k} \to C(\overline{G}_{\reg};F)^{\overline{G}\times\mathcal{G}_k}\]
as follows:
\[\Ind_{\overline{H}}^{\overline{G}}(f)(g) = \frac{1}{|\overline{H}|}\sum_{x \in \overline{G}}f(xgx^{-1})\]
for any $f \in C(\overline{H}_{\reg};F)^{\overline{H}\times\mathcal{G}_k}$ and $g \in \overline{G}_{\reg}$, with the understanding that $f(u) = 0$ if $u \notin \overline{H}$.
Consider the following diagram:
\[
\xymatrix{
FK_0(\mathcal{F}_d(G)) \ar@{=}[r] & F\mathcal{G}_0(kG) \ar[r]^{\theta_N} & F\mathcal{G}_0(k\overline{G})\ar[r]^{\varphi} & C(\overline{G}_{\reg};F)^{\overline{G}\times\mathcal{G}_k} \\
FK_0(\mathcal{F}_e(H)) \ar[u]^{\Ind_H^G}\ar@{=}[r] & F\mathcal{G}_0(kH) \ar[u]^{\Ind_H^G}\ar[r]_{\theta_{H\cap N}} & F\mathcal{G}_0(k\overline{H}) \ar[u]_{\Ind_{\overline{H}}^{\overline{G}}} \ar[r]_{\varphi} & C(\overline{H}_{\reg};F)^{\overline{H}\times\mathcal{G}_k} \ar[u]_{\Ind_{\overline{H}}^{\overline{G}}}.
}
\]
The middle square commutes by functoriality of $\mathcal{G}_0$ and the right-hand square commutes by the well-known formula \cite[Theorem 12 and Exercise 18.2]{Serre2} for the character of an induced representation. So the whole diagram commutes.

Finally, using the isomorphisms $\pi_N^\ast$ and $\pi_{H\cap N}^\ast$ which feature in (\ref{PfBmWt}) we can define the required map $\iota_N$ to be the map which makes the following diagram commute:
\[
\xymatrix{
C(\overline{G}_{\reg};F)^{\overline{G}\times\mathcal{G}_k} \ar[r]^{\pi_N^\ast} & C(G_{\reg};F)^{G\times\mathcal{G}_k} \\
C(\overline{H}_{\reg};F)^{\overline{H}\times\mathcal{G}_k} \ar[r]_{\pi_{H\cap N}^\ast}\ar[u]^{\Ind_{\overline{H}}^{\overline{G}}} & C(H_{\reg};F)^{H\times\mathcal{G}_k} \ar[u]_{\iota_N}.
}
\]
The commutative diagram appearing in (\ref{PfBmWt}) now shows that
\[\pi_N^\ast\circ \varphi = \varphi\circ\lambda_N\quad \mbox{and}\quad \pi_{H\cap N}^\ast\circ \varphi = \varphi\circ\lambda_{H\cap N}\]
and the result follows by pasting the above diagrams together.
\end{proof}

\subsection{The case when $G$ is virtually abelian}
\label{VirtAb}
We can now apply the theory developed above. If $g\in G_{\reg}$ let $\delta_g : G_{\reg} \to F$ be the function which takes the value $1$ on the $G\times \mathcal{G}_k$-orbit of $g$ and is zero elsewhere. Note that $\delta_g$ is locally constant by Corollary \ref{HallBij}.

\begin{thm} Suppose that $G$ is virtually abelian. Then the upper bound of Theorem \ref{UpBnd} is always attained:
\[\rk \alpha_i = |(G\times \mathcal{G}_k)\backslash S_i|\]
for all $i=0,\ldots,d$.
\end{thm}
\begin{proof}
Fix the integer $i$ and fix $g\in G_{\reg}$ such that $\dim C_G(g) \leq i$. As $\rk \alpha_i = \rk (\rho_N \circ \alpha_i)$ by Corollary \ref{kGbar}, it will be sufficient to show that
\[\delta_g \in \im (\rho_N \circ \alpha_i).\]

Because we are assuming that $G$ is virtually abelian, the uniform pro-$p$ subgroup $N$ is abelian. Thus $N$ is a free $\Zp$-module of rank $d = \dim G$. By considering the conjugation action of $g$ on $N$ we see that the submodule of fixed points
\[C_N(g) = \{x \in N : gx = xg\} = N \cap C_G(g)\]
has a unique $g$-invariant $\Zp$-module complement in $N$ which we will denote by $L$. Note that $L$ is a \emph{subgroup} of $G$ as $N$ is abelian. We could alternatively have defined $L$ as the isolator of $[N,\langle g\rangle]$ in $N$.

Let $H$ be the closed subgroup of $G$ generated by $L$ and $g$. Because $g$ normalises $L$, $H$ is isomorphic to a semi-direct product of $L$ with the finite group $\langle g \rangle$:
\[H = L \rtimes \langle g \rangle.\]

Let $\epsilon_g : H_{\reg} \to F$ be the locally constant function which is $1$ on the $H\times\mathcal{G}_k$-orbit of $g$ inside $H_{\reg}$ and zero elsewhere. Note that $\epsilon_g$ is constant on the cosets of the open uniform subgroup $H\cap N = L$ of $H$.

By construction, $g$ acts without nontrivial fixed points on $L$ by conjugation. So if $\beta_g$ denotes this action, then
\[D_g := \det(1 - \beta_g)\in F\]
is a nonzero constant. In view of the commutative diagram which appeared in the proof of Proposition \ref{rkalpha0} the element $X:=\varphi^{-1}(\epsilon_g) \in FK_0(\mathcal{F}_0(H))$ satisfies
\[(\rho_{H\cap N}\circ\alpha_0)(X) = D_g\epsilon_g.\]

Using the definition of the map $\iota_N : C(H_{\reg};F)^{H\times\mathcal{G}_k} \to C(G_{\reg};F)^{G\times\mathcal{G}_k}$ of Proposition \ref{IndComDiag} we may calculate that there exists a nonzero constant $A_g\in F$ such that for all $y\in G_{\reg}$ we have
\[\iota_N(\epsilon_g)(y) = A_g \cdot \delta_g(y). \]
We can now apply Proposition \ref{IndComDiag} and obtain
\[(\rho_N\circ \alpha_{d-e})\left(\Ind_H^GX\right) = \iota_N(\rho_{H\cap N}(\alpha_0(X)))=\iota_N\left(D_g\epsilon_g\right) = D_g A_g \delta_g.\]
But $d-e = \dim G - \dim H = \dim C_G(g) \leq i$ and $D_gA_g \neq 0$ so
\[\delta_g \in \im (\rho_N\circ\alpha_{d-e}) \subseteq \im (\rho_N \circ \alpha_i)\]
as required.
\end{proof}

\section{An example}

\subsection{Central torsion modules}
\label{CentTors}
Suppose that $z$ is a central element of $G$ contained in $N$. Write $Z$ for the closed central subgroup of $G$ generated by $z$.

Because $N$ is torsion-free by \cite[Theorem 4.5]{DDMS}, $Z$ is
isomorphic to $\Zp$. Hence $kZ$ is isomorphic to the power series
ring $k[[z-1]]$ and its maximal ideal is generated by $z-1$. Because
$kG$ is a flat $kZ$-module, it follows that $z-1$ generates the
kernel of the map $kG \twoheadrightarrow k[[G/Z]]$ as an ideal in
$kG$.

\begin{lem} Let $M \in \mathcal{M}(kG)$ and suppose that $M.(z-1) = 0$ so that $M$ is also a right $k[[G/Z]]$-module. Then as right $k[[G/Z]]$-modules, we have
\[\Tor_n^{kG}(M, k[[G/Z]]) \cong \left\{ \begin{array}{cc}M & \mbox{ if } \quad n=0 \mbox{ or }n=1 \\ 0 & \mbox{ otherwise.} \end{array} \right.
\]
\end{lem}
\begin{proof}
Because $z-1$ is not a zero-divisor in $kG$,
\[0 \to kG \stackrel{z-1}{\longrightarrow} kG \to k[[G/Z]] \to 0\]
is a free resolution of $k[[G/Z]]$ as a left $kG$-module. Hence the complex
\[0 \to M \stackrel{z-1}{\longrightarrow} M \to 0 \]
computes the required modules $\Tor_n^{kG}(M,k[[G/Z]])$. The result follows because $M.(z-1)=0$.
\end{proof}

\subsection{Proposition}
\label{ThetaZero}
Let $M\in\mathcal{M}(kG)$ and suppose that $M.(z-1)^a = 0$ for some $a\geq 1$. Then $\theta_N[M] = 0$.
\begin{proof}
Suppose first that $a=1$ so that $M$ is killed by $z-1$. We recall the base-change spectral sequence for $\Tor$ \cite[Theorem 5.6.6]{Wei} associated with a ring map $f : R \to S$:
\[E^2_{ij}=\Tor_i^S(\Tor_j^R(A,S),B) \Longrightarrow \Tor_{i+j}^R(A,B).\]
This is a first quadrant convergent homological spectral sequence. We apply it to the map $R := kG \twoheadrightarrow k[[G/Z]] =: S$ with $A := M$ and $B := k\overline{G} = k[G/N]$ --- note that $B$ is a left $k[[G/Z]]$-module because we are assuming that $Z \leq N$.

By Lemma \ref{CentTors} this spectral sequence is concentrated in rows $j=0$ and $j=1$, so we may apply \cite[Exercise 5.2.2]{Wei} and obtain a long exact sequence
\[\begin{array}{rcccccl} \cdots \to &\Tor_{n+1}^{kG}(M,k\overline{G}) &\to &\Tor_{n+1}^{k[[G/Z]]}(M,k\overline{G}) &\to &\Tor_{n-1}^{k[[G/Z]]}(M,k\overline{G})&\to\\

\to &\Tor_n^{kG}(M,k\overline{G}) &\to &\Tor_n^{k[[G/Z]]}(M,k\overline{G}) &\to &\Tor_{n-2}^{k[[G/Z]]}(M,k\overline{G})&\to\cdots.
\end{array}
\]
We can now apply Lemma \ref{HomEul} and deduce that
\[\theta_N[M] = \sum_{n=0}^d (-1)^n [\Tor_n^{kG}(M,k\overline{G})] = 0 \in \mathcal{G}_0(k\overline{G})\]
as required. The general case follows quickly by an induction on $a$.
\end{proof}

\begin{cor} Let $M$ be as above. Then its Euler characteristic is trivial:
\[\chi(G,M) = 1.\]
\end{cor}
\begin{proof} This follows from Definition \ref{GenEulCh}.
\end{proof}

Using this result we now show that the upper bound for $\rk \alpha_i$ given in Theorem \ref{UpBnd} is not \emph{always} attained.

\subsection{Example}
\label{CEx}
Let $p$ be odd and let $N$ be a \emph{clean Heisenberg pro-$p$ group} \cite[Definition 4.2]{W} of dimension $2r+1$. By definition, $N$ has a topological generating set $\{x_1,\ldots,x_r,y_1,\ldots,y_r,z\}$ and relations $[x_i,y_i]=z^p$ for each $i=1,\ldots, r$, all other commutators being trivial. Note that $N$ is uniform.

The presentation for $N$ given above makes it possible to define an automorphism $\gamma$ of $N$ which fixes $z$ and sends all the other generators to their inverses:
\[\gamma(x_i) = x_i^{-1}, \gamma(y_i) = y_i^{-1}, \gamma(z) = z.\]

Now let $G$ be the semidirect product of $N$ with a cyclic group $\langle g\rangle $ of order $2$, where the conjugation action of $g$ on $N$ is given by the automorphism $\gamma$. Thus the $p'$-part of $|G/N|$ is equal to $2$ so the Galois group $\mathcal{G}_k$ defined in (\ref{preg}) is trivial for any finite field $k$ of characteristic $p$.

Now $(G/N)_{\reg} = G/N$ has two $G/N$-conjugacy classes, so by Proposition \ref{HallBij} there are just two $G\times\mathcal{G}_k$-orbits on $G_{\reg}$, represented by the elements $1$ and $g$.

By construction, the endomorphism $\gamma = \Ad(g)$ has eigenvalue $-1$ with multiplicity $2r$ and eigenvalue $1$ with multiplicity $1$:
\[\Psi(g) = (1+t)^{2r}(1-t),\]
so $\dim C_G(g) = 1$ and of course $\dim C_G(1) = 2r + 1$. Hence in the notation of (\ref{UpBndPf})
\[|(G\times\mathcal{G}_k)\backslash S_i| = \left\{ \begin{array}{cc}0 & \mbox{ if }  i=0, \\  1 &\mbox{if  } 1\leq i \leq 2r, \\ 2 & \mbox{if }i=2r+1. \end{array} \right.\]

Now if $M$ is a finitely generated $kG$-module with $d(M) \leq r$ then $M$ is killed by some power of $z-1$ by \cite[Theorem B]{W}. It follows from Proposition \ref{ThetaZero} that $\theta_N[M] = 0$, so $\rk \alpha_i = 0$ for all $i\leq r$ --- thus the upper bound given in Theorem \ref{UpBnd} is not attained for all values of $i$ between $i=1$ and $i=r$.

\end{document}